\documentclass[aop,MSNbibl,seceqn,citesort,dvips]{arximspdf}
\usepackage{subenv}

\doi{10.1214/11-AOP662}
\volume{40}
\issue{4}
\pubyear{2012}
\firstpage{1675}
\lastpage{1714}

\makeatletter

\newcommand{\iint}{\int\!\!\int}

\newcommand{\Overline}{\bar}

\newcommand{\Si}{\operatorname{Si}}

\newcommand{\R}{\mathbf{R}}
\newcommand{\Z}{\mathbf{Z}}
\newcommand{\C}{\mathbf{C}}

\newtheorem{theorem}{Theorem}[section]
\newtheorem{corollary}[theorem]{Corollary}
\newtheorem{lemma}[theorem]{Lemma}
\newtheorem{proposition}[theorem]{Proposition}

\newproclaim{assumption}[theorem]{Assumption}
\newproclaim{remark}[theorem]{Remark}

\newcommand{\E}{{\mathbb E}}

\newcommand{\ot}{\otimes}
\newcommand{\ta}{{\widetilde{\alpha}}}
\newcommand{\hD}{\widehat{D}}

\newcommand{\BV}{\mathrm{BV}}
\newcommand{\tF}{\Overline{F}}

\newcommand{\psio}{\psi}
\newcommand{\psiog}{\psi_\gamma}
\newcommand{\psioug}{\psi^\gamma}
\newcommand{\psiouz}{\psi^\zeta}
\newcommand{\psiogx}{\psi_\gamma^\chi}

\newcommand{\psit}{\widetilde{\psi}}
\newcommand{\psitg}{\widetilde{\psi}_\gamma}
\newcommand{\psitz}{\widetilde{\psi}_\zeta}
\newcommand{\psitlx}{\widetilde{\psi}_\chi}
\newcommand{\psitug}{\widetilde{\psi}^\gamma}
\newcommand{\psituz}{\widetilde{\psi}^\zeta}
\newcommand{\psitgx}{\widetilde{\psi}_\gamma^\chi}

\newcommand{\vb}{\bar v}
\newcommand{\vbg}{\bar v^\gamma}
\newcommand{\vt}{\widetilde v}
\newcommand{\vtg}{\widetilde v^\gamma}

\newcommand{\thk}{\theta_k}
\newcommand{\thtk}{\widetilde{\theta}_k}
\newcommand{\thtmk}{\widetilde{\theta}_{-k}}

\newcommand{\one}{\mathbf{1}}

\newcommand{\tr}{\operatorname{tr}}

\newcommand{\CE}{\mathcal{E}}
\newcommand{\CR}{\mathcal{R}}
\newcommand{\CH}{\mathcal{H}}
\newcommand{\CC}{\mathcal{C}}

\makeatother

\begin{document}
\begin{frontmatter}

\title{A spatial version of the It\^{o}--Stratonovich correction}
\runtitle{A spatial version of the It\^{o}--Stratonovich correction}

\begin{aug}
\author[A]{\fnms{Martin} \snm{Hairer}\thanksref{t1}\ead[label=e1]{M.Hairer@Warwick.ac.uk}} and
\author[B]{\fnms{Jan} \snm{Maas}\corref{}\thanksref{t2}\ead[label=e2]{maas@iam.uni-bonn.de}}
\runauthor{M. Hairer and J. Maas}
\affiliation{University of Warwick and University of Bonn}
\address[A]{Department of Mathematics\\
University of Warwick\\
Coventry CV4 7AL\\
United Kingdom\\
\printead{e1}} 
\address[B]{Institute for Applied Mathematics\\
University of Bonn\\
Endenicher Allee 60\\
53115 Bonn\\
Germany\\
\printead{e2}}
\end{aug}

\thankstext{t1}{Supported by the EPSRC Grants EP/E002269/1 and EP/D071593/1,
a Wolfson Research Merit Award of the Royal Society and a Philip Leverhulme prize
of the Leverhulme Trust.}

\thankstext{t2}{Supported by Rubicon Grant 680-50-0901
of the Netherlands Organisation for Scientific Research (NWO).}

\received{\smonth{11} \syear{2010}}
\revised{\smonth{3} \syear{2011}}

%
\begin{abstract}
We consider a~class of stochastic PDEs of Burgers type in spatial
dimension $1$,
driven by space--time white noise. Even though it is well known that
these equations
are well posed, it turns out that if one performs a~spatial
discretization of the nonlinearity
in the ``wrong'' way, then the sequence of approximate equations does
converge to a~limit,
but this limit exhibits an additional correction term.

This correction term is proportional to the local quadratic
cross-variation (in space) of
the gradient of the conserved quantity with the solution itself. This
can be understood as a~consequence of the fact that
for any fixed time, the law of the solution is locally equivalent to
Wiener measure, where
space plays the role of time. In this sense, the correction term is
similar to the usual
It\^{o}--Stratonovich correction term that arises when one considers
different temporal
discretizations of stochastic ODEs.
\end{abstract}

%

\pdfkeywords{60H15, 35K55, 60H30, 60H35, Ito--Stratonovich correction,
stochastic Burgers equation, spatial discretizations,
Wiener chaos}

\end{frontmatter}

\section{Introduction}

In this work, we give a~rigorous analysis of the behavior of
stochastic Burgers equations
in one spatial dimension under various approximation schemes. It was
recently argued in \cite{Jochen} that
if the approximation scheme fails to satisfy a~certain symmetry
condition, then
one expects the approximations to converge to a~modified equation, with
the appearance of
an additional correction term in the limit. This correction term is
somewhat similar to the It\^{o}--Stratonovich correction
that appears in the study of SDEs when one compares centred and
one-sided approximations.
The present article provides a~rigorous justification of the results
observed in \cite{Jochen}, at least in the
case where the nonlinearity of the equation is of gradient type, and
therefore the limiting equation is well posed
in a~classical sense.

More precisely, we will consider in this work stochastic PDEs of the
form
%
%
\begin{equation}\label{eSPDE}
\partial_t u(x,t)
= \nu\,\partial_x^2 u(x,t) + F(u(x,t)) + (\nabla G)( u(x,t) )\,
\partial_x u(x,t) + \xi(x,t) ,\hspace*{-32pt}
\end{equation}
where $u = u(x,t)$ is an $\R^n$-valued function, with $x \in[0, 2\pi
]$, $t \geq0$. In this equation, $\nu> 0$ is a~positive constant, the
functions $F, G \dvtx\R^n \to\R^n$ are assumed to be $\CC^3$, and the
stochastic forcing term $\xi$ consists of independent space--time white
noises in each component of $\R^n$. For the sake of simplicity, we
endow this equation with periodic boundary conditions, but we do not
expect this
to alter our results significantly.

Endowed with an initial condition $u_0 \in\CC([0,2\pi];\R^n)$,
(\ref{eSPDE}) is locally well posed \cite{Gyo98}, provided that we
rewrite the nonlinearity as $\partial_x G(u)$ and consider solutions
either in the weak or the mild form \cite{ZDP}. (Note that our noise
term is \textit{not} the gradient of space--time white noise, as
in \cite{BertGia97}. Therefore, our solutions are actually $\alpha
$-H\"{o}lder continuous functions for all $\alpha< \frac12$.)
The aim of this article is to show that this well-posedness is much
less stable than
one may imagine at first. Indeed, if we set
\[
D_\varepsilon^+ u(x) = {u(x+\varepsilon) - u(x) \over\varepsilon} ,
\]
and consider the family $u_\varepsilon$ of solutions to the
approximate equation
\[
\partial_t u_\varepsilon= \nu\,\partial_x^2 u_\varepsilon+
F(u_\varepsilon) + \nabla G(u_\varepsilon) D_\varepsilon^+
u_\varepsilon+ \xi,
\]
then our main result, Theorem~\ref{thmmain-one} below, implies that
$u_\varepsilon\Rightarrow\bar u$, where $\bar u$ is the solution
to (\ref{eSPDE}), but with $F$ replaced by
%
%
\begin{equation}\label{ecorF}
\bar F(u) = F(u) - {1\over4\nu} \Delta G(u) ,
\end{equation}
where $\Delta$ is the usual Laplacian on $\R^n$.
\begin{remark}
The correction term in (\ref{ecorF}) is reminiscent of the Wong--Zakai
correction \cite{WZ}, which arises if the driving Brownian motion in a~stochastic ODE or PDE is approximated by stochastic processes of
bounded variation. This correction term is due to the \textit{temporal}
roughness of the driving Brownian motion and does not appear if the
noise is additive.

Our correction term is a~consequence of the \textit{spatial} roughness
of the solutions and appears even though we consider SPDEs with
additive noise.
In fact, an explicit calculation allows to check that the local
quadratic variation (in space) of the solution $u$ to (\ref{eSPDE})
is precisely given by $1/(2\nu)$. Therefore, one can interpret the
correction term appearing in (\ref{ecorF}) as precisely being equal
to $-{1\over2}$ times the quadratic covariation between $u$ and
$\nabla G(u)$. Recall that this is exactly the correction term that appears
when one switches between It\^{o} and Stratonovich integral in the
usual setting of stochastic~calculus.
See also \cite{Jochen} for a~heuristic argument for computing the correction~term.
\end{remark}
\begin{remark}
This correction term is a~purely stochastic effect and is completely
unrelated to the fact that our
discretization scheme is not an upwind scheme (see \cite
{Upwind,BookNA}). In the absence of noise, we would still
have the regularizing property from the nonvanishing viscosity, so that
pretty much any ``reasonable''
numerical scheme would converge to the correct solution.
\end{remark}

If $D_\varepsilon^+$ is replaced by $D_\varepsilon^-$, defined
by $D_\varepsilon^- u(x) = (u(x) - u(x-\varepsilon)) /
\varepsilon$, then a~similar result is true, but the sign in front of the
correction term in (\ref{ecorF}) changes.
We will actually consider a~much more general class of approximations
to (\ref{eSPDE}), where we also allow
both the linear operator $\partial_x^2$ and the noise term $\xi$ to be
replaced by approximate versions that are still translation-invariant,
but modified at the lengthscale $\varepsilon$.

\subsection{Statement of the main result} \label{secmain}

For $\varepsilon> 0$, we consider approximating stochastic PDEs of the type
\[
\partial_t u_\varepsilon= \nu\Delta_\varepsilon u_\varepsilon+
F(u_\varepsilon)
+ \nabla G(u_\varepsilon) D_\varepsilon u_\varepsilon+ \xi
_\varepsilon.
\]
Since our system is invariant under spatial translations, it seems
natural to restrict ourselves to a~class of approximations that enjoys the same property. Throughout this
article, we will therefore
use approximate differential operators $\Delta_\varepsilon$
and~$D_\varepsilon$, as well as an approximate space--time white noise
$\xi_\varepsilon$ given by their Fourier transforms:
\begin{eqnarray*}
\widehat{\Delta_\varepsilon u}(k) &=& - k^2 f(\varepsilon|k|)\widehat
u(k),\qquad
\widehat{D_\varepsilon u}(k) = ik g(\varepsilon k)\widehat u(k),\\
\widehat{\xi}_\varepsilon(k) &=& h(\varepsilon|k|) \widehat
\xi(k).
\end{eqnarray*}
Several natural discretizations arising in numerical analysis are of
this form (see the examples below Theorem~\ref{thmmain-one}).
We will make the following standing assumptions on these objects.
\begin{assumption}
The function $f \dvtx[0,\infty) \to[0,+\infty]$ is twice differentiable
at $0$ with $f(0) = 1$ and $f'(0) = 0$.
Furthermore, there exists $q\in(0,1]$ such that $f(k) \ge q$ for all
$k > 0$.
\end{assumption}

If $f(k) = +\infty$ for some values of $k$, we use the convention
$\exp(-t \infty) = 0$ for every $t>0$. In this case,
the semigroup generated by $\Delta_\varepsilon$ is not strongly
continuous, but this is of no consequence for our analysis.
\begin{assumption}
There exists a~signed Borel measure $\mu$ such that
\[
\int_\R e^{ikx} \mu(dx) = ik g(k) ,
\]
and such that
%
%
\begin{equation} \label{eqmoments}\qquad
\mu(\R) = 0 ,\qquad
|\mu|(\R) < \infty,
\int_{\R} x \mu(dx) = 1 ,\qquad
\int_{\R} |x|^4 |\mu|(dx) < \infty.
\end{equation}
In particular, we have $ (D_\varepsilon u)(x) := \frac{1}{\varepsilon
} \int_\R u(x + \varepsilon y) \mu(dy)$.\vadjust{\goodbreak}
\end{assumption}
\begin{assumption}
The function $h$ is bounded and such that $h^2/f$ is of bounded variation.
Furthermore, $h$ is twice differentiable at the origin with \mbox{$h(0)
= 1$}
and $h'(0) = 0$.
\end{assumption}

Let $\bar u$ be the solution to the equation
%
%
\begin{equation}\label{eSPDE2}
\partial_t \bar u = \nu\,\partial_x^2 \bar u + \tF(\bar u)
+ \nabla G(\bar u) \,\partial_x \bar u + \xi.
\end{equation}
In this equation,
\[
\tF:= (F - \Lambda\Delta G)
\]
and $\Lambda\in\R$ is a~correction constant given by
%
%
\begin{equation}\label{edefLambda}
\Lambda\stackrel{\mathrm{def}}{=}\frac{1}{2\pi\nu} \int_{\R_+}
\int_\R\frac{(1- \cos(yt))h^2(t)}{t^2 f(t)} \mu(dy) \,dt .
\end{equation}
Note that a~straightforward calculation shows that $\Lambda$ is indeed
well defined, as a~consequence of the fact that
$h^2 / f$ is bounded by assumption and that~$|\mu|$ has a~finite
second moment.

Before we state our main result, note that the equation (\ref{eSPDE2})
is locally well posed in $L^\infty$,
see \cite{BCF91,BCJ94,DDT94,Gyo98,Hai09}.
As a~consequence, it has a~well-defined blow-up time $\tau_*$
(possibly infinite) such that, almost surely,
${\lim_{t \to\tau_*}}\|\bar u(t)\|_{L^\infty} = +\infty$ on the
event $\{\tau_* < \infty\}$.
With this notation, we are now ready to state the main result of this paper.
\begin{theorem}\label{thmmain-one}
Let $\kappa> 0$ and let $u_\varepsilon$ and $\bar u$ have initial
conditions as in Theorem~\ref{thmmain}. There exists a~sequence of
stopping times $\tau_\varepsilon$ satisfying $\lim_{\varepsilon\to
0}\tau_\varepsilon= \tau_*$
in probability, and such that
\[
\lim_{\varepsilon\to0} {\mathbb P}\Bigl( \sup_{t \leq\tau
_\varepsilon} \| u_\varepsilon(t) - \bar u(t)\|_{L^\infty} >
\varepsilon^{1/8 - \kappa} \Bigr) = 0 .
\]
\end{theorem}
\begin{remark}
In order to avoid further technical complications, we consider
sequences of initial conditions that have the property that the initial
condition for $u_\varepsilon$ ``behaves like'' the solution
$u_\varepsilon(t)$ for positive times. In fact, the initial condition
for $u_\varepsilon$ is a~smooth perturbation of the stationary
solution to the linearized equation for $u_\varepsilon$. We refer to
Section~\ref{secproof} for more details.
\end{remark}

Before we proceed, we list some of the most common examples of
discretizations that do fit our framework.
For $a,b \ge0$ with $a+b>0$, it is natural to discretize the
derivative operator by
choosing
\[
\mu:= \frac{\delta_{a} - \delta_{-b}}{a + b} .
\]
This is also the discretization that was considered in \cite{Jochen}.
As far as the discretizations
of the noise and the Laplacian are concerned, there are at least three
natural choices.

\subsubsection*{No discretization} This is the case $f=h=1$ where
only the nonlinearity is discretized.
With this choice, one can check that the correction factor is given by
$\Lambda= \frac{1}{4\nu}\frac{a-b}{a+b}$.

\subsubsection*{Finite difference discretization} In this case, we
divide the interval $[0,2\pi]$ into $N$
equally sized intervals. For convenience, we assume that $N$ is odd and
we set
\[
(\Delta_\varepsilon u)(x) = \frac{1}{\varepsilon^2}
\bigl(u(x+\varepsilon) + u(x-\varepsilon) -2u(x)\bigr) ,\qquad
\varepsilon= {2\pi\over N} .
\]
We furthermore identify a~function $u$ with the trigonometric
polynomial of degree $(N-1)/2$ agreeing with $u$ at the
gridpoints.
This corresponds to the choice
\[
f(k) = \cases{
\displaystyle\frac4{k^2} \sin^2({k/2}), &\quad$k \in[0,\pi)$,\vspace
*{2pt}\cr
+ \infty, &\quad$k \in[\pi, \infty)$,}\qquad
h = \one_{[0,\pi)}.
\]
The natural choice for the discretization of the derivative operator in
this case is to choose $a$ and $b$
to be integers, so that discretization takes place on the gridpoints.
With this choice, it can be shown that
the correction factor is identical to that obtained in the previous
case. Note however that this is \textit{not} the case
if the discretization of the derivative operator is not adapted to the gridsize.

\subsubsection*{Galerkin discretization} In this case, we approximate
$\Delta$ and $\xi$ by only keeping
those Fourier modes that appear in the approximation by trigonometric
polynomials. This corresponds to the choice
\[
f(k) = \cases{
1, &\quad$k \in[0,\pi)$,\cr
+ \infty, &\quad$k\in[\pi, \infty)$,}\qquad
h = \one_{[0,\pi)} .
\]
The correction factor $\Lambda$ is then given by
\[
\Lambda= {\cos(\pi a) + \pi a \Si(\pi a) - \cos(\pi b) - \pi b\Si
(\pi b) \over2\pi^2 \nu(a+b)} ,
\]
where $\Si t = \int_0^t {\sin x \over x} \, dx$.\vspace*{1pt}

The rest of this paper is structured as follows. In Section \ref
{secproof}, we introduce notation,
we give a~refined formulation of the main result and present an outline
of the proof of the main
result (Theorem~\ref{thmmain}). In Section~\ref{secingredients}, we
prove several useful bounds
on the approximating semigroups and the approximations of the gradient.
Section~\ref{secprobtools}
is devoted to several estimates for stochastic convolutions, the most
crucial one
being Proposition~\ref{propsecond-chaos}, which is responsible for
the correction term appearing in the
limiting equation. Most of the work is performed in Section \ref
{secapprox}, where
convergence of various approximating equations is proved.

\section{Proof of the main result}
\label{secproof}

In order to shorten notation, we introduce the semigroups $S$ and
$S_\varepsilon$, defined as rescaled versions of the
heat semigroup and its approximation:
\[
S(t) \stackrel{\mathrm{def}}{=}e^{-t(1 - \nu\partial_x^2)} ,\qquad
S_\varepsilon(t) \stackrel{\mathrm{def}}{=}e^{-t(1 - \nu\Delta
_\varepsilon
)} ,
\]
where we define $S_\varepsilon$ by Fourier analysis, that is,
\[
\widehat{S_\varepsilon u}(k) = e^{-t (1 + \nu k^2 f(\varepsilon
|k|))}\widehat u(k) ,
\]
making use of the convention $e^{-\infty} = 0$.

Since we will always work with the mild formulation, it will be
convenient to have a~notation for the convolution
(in time) of a~function with one of the semigroups. We will henceforth write
\[
(S*w)(t) \stackrel{\mathrm{def}}{=}\int_0^t S(t-s) w(s) \,ds .
\]
Let $(W(t))_{t \in\R}$ be a~two-sided cylindrical Wiener process on
$\CH\stackrel{def}{=}L^2([0,2\pi],\allowbreak\R^n)$ (see \cite{ZDP,Hai09} for precise
definitions)
and let $Q_\varepsilon$ be the bounded operator on~$\CH$ defined
as a~Fourier multiplier by
\[
\widehat{Q_\varepsilon u}(k) = h(\varepsilon|k|) \widehat{u}(k) .
\]
(We assume\vspace*{2pt} that it acts independently on each component.)
Finally, we define the $\CH$-valued processes $\psio$ and $\psit$ by
\[
\psio(t) = \int_{-\infty}^t S(t-s) \,dW(s) ,\qquad
\psit(t) = \int_{-\infty}^t S_\varepsilon(t-s) Q_\varepsilon
\,dW(s) ,
\]
so that, in the notation of the previous section, they are the
\textit{stationary} solutions
to the linear equations
\[
\partial_t \psio= (\nu\partial_x^2 - 1) \psio+ \xi,\qquad
\partial_t \psit= (\nu\Delta_\varepsilon- 1) \psit+ \xi
_\varepsilon.
\]
With this notation at hand, we can rewrite the equations for
$\bar u$ and $u_\varepsilon$ in the mild form as
%
%
\begin{eqnarray} \label{equ-0}
\bar u(t) & = & S(t) v_0 + \psio(t) + S * \bigl( \tF(\bar u) + \nabla
G(\bar u) \,\partial_x \bar u
\bigr)(t) ,\\ \label{equ-eps}
u_\varepsilon(t) & = & S_\varepsilon(t) v_0 + \psit(t)
+ S_\varepsilon* \bigl( F(u_\varepsilon) + \nabla G(u_\varepsilon)
D_\varepsilon u_\varepsilon
\bigr)(t) .
\end{eqnarray}

\begin{remark}
Note that we have used here a~common initial condition~$v_0$ for the
difference $\bar u - \psio$ and $u_\varepsilon- \psit$.
As a~consequence, the two equations \textit{do not} start with the
same initial condition! However, as $\varepsilon\to0$, the initial
condition of $u_\varepsilon$ converges to that of $\bar u$. The reason
for not starting with the same initial condition is mostly
of technical nature.
\end{remark}

It will be convenient to define for any $0 < \gamma< \chi$,
\[
\psiog\stackrel{\mathrm{def}}{=}(I - \Pi_{\varepsilon^{-\gamma}})
\psio,\qquad
\psioug\stackrel{\mathrm{def}}{=}\Pi_{\varepsilon^{-\gamma}}
\psio,\qquad
\psiogx\stackrel{\mathrm{def}}{=}(\Pi_{\varepsilon^{-\chi}} - \Pi
_{\varepsilon
^{-\gamma}}) \psio.\vadjust{\goodbreak}
\]
The expressions $\psitg$, $\psitug$ and $\psitgx$ are defined analogously.
Here $\Pi_N$ denotes the projection onto the low-frequency components
of the Fourier expansion, defined by $\Pi_N e_n \stackrel{\mathrm
{def}}{=}\one
_{|n| \leq
N} e_n$, where $e_n(x) = (2\pi)^{-{1/2}} e^{inx}$.

We set
\[
\vb\stackrel{\mathrm{def}}{=}\bar u - \psio,\qquad
\vt\stackrel{\mathrm{def}}{=}u_\varepsilon- \psit.
\]
In the proof, it will be convenient to work with the functions $\vbg$
and $\vtg$ defined by
\[
\vbg\stackrel{\mathrm{def}}{=}\vb+ \psioug= \bar u - \psiog,\qquad
\vtg\stackrel{\mathrm{def}}{=}\vt+ \psitug= u_\varepsilon- \psitg.
\]
It follows from (\ref{equ-0}) and (\ref{equ-eps}) that these
functions satisfy the following equations:
%
%
\begin{eqnarray}\qquad
\label{eqv-0}
\vbg(t)& = & S(t) v_0
+ \psioug(t)
+ S * \bigl( \tF( \vbg+ \psiog)
+ \partial_x \bigl(G( \vbg+ \psiog
)\bigr) \bigr)(t) ,\\
\label{eqv-eps}
\vtg(t) & = & S_\varepsilon(t) v_0
+ \psitug(t)\nonumber\\[-8pt]\\[-8pt]
&&{} + S_\varepsilon* \bigl(
F(\vtg+ \psitg)
+ \nabla G(\vtg+ \psitg)
D_\varepsilon(\vtg+ \psitg)
\bigr)(t).\nonumber
\end{eqnarray}

For large parts of this article, it will be convenient to work in the
fractional Sobolev space
$H^\alpha$ for some $\alpha> {1\over2}$,\vspace*{1pt} so that $H^\alpha
\subset
L^\infty$. Recall that $H^\alpha$ denotes the space of (equivalence
classes of) functions $u = \sum_{j \in\Z} u_k e_k$ on $[0,2\pi]$
with $u_k \in\C$, for which
\[
\| u \|_\alpha^2 := \sum_{k \in\Z} |u_k|^2 (1 + k^2)^\alpha
< \infty.
\]
Furthermore, we will need
to use a~high-frequency cut-off, which will smoothen out the solutions
at a~scale $\varepsilon^{\chi}$ for some $\chi> 1$.
It turns out that a~reasonable choice for these parameters is given by
%
%
\begin{equation}\label{eexponents}
\alpha= \tfrac34 , \qquad\gamma= \tfrac13 ,\qquad
\chi= \tfrac32 ,
\end{equation}
and we will fix these values from now on.
With this notation at hand, the following theorem, which is
essentially a~more precise reformulation of Theorem~\ref{thmmain-one},
is a~more precise statement of our main result. Here and in the rest of
the paper we write $\|u\|_{\beta}$ to denote the norm of an element
$u$ in the fractional Sobolev space $H^\beta$ for $\beta\in\R$.
\begin{theorem}\label{thmmain}
Let $\kappa> 0$ be an arbitrary (small) exponent and let \mbox{$v_0 \in
H^{\beta}$} for all $\beta< \frac32$.
There exists a~sequence of stopping times $\tau_\varepsilon$
satisfying $\tau_\varepsilon\to\tau_*$ in probability as
$\varepsilon\to0$, such that
%
%
\begin{equation}\label{ebounduu}
\lim_{\varepsilon\to0} {\mathbb P}\Bigl( \sup_{t \leq\tau
_\varepsilon} \| u_\varepsilon(t) - \bar u(t)\|_{L^\infty} >
\varepsilon^{1/8 - \kappa} \Bigr) = 0 .
\end{equation}
In fact, we have the bounds
%
%
\begin{eqnarray}
\label{eboundvv}
\lim_{\varepsilon\to0} {\mathbb P}\Bigl( \sup_{t \leq\tau
_\varepsilon} \| \vtg(t) - \vbg(t) \|_{\alpha} > \varepsilon
^{1/8 - \kappa} \Bigr) & = & 0 ,\\
\label{eboundpsipsi}
\lim_{\varepsilon\to0} {\mathbb P}\Bigl( \sup_{t \leq\tau
_\varepsilon} \| \psitg(t) - \psiog(t) \|_{L^\infty} > \varepsilon
^{1/2 - \kappa} \Bigr) & = & 0 .
\end{eqnarray}
\end{theorem}
\begin{remark}\label{reminitial-value}
We emphasize again that the initial conditions $\bar u(0)$ and
$u_\varepsilon(0)$ are slightly different. In fact, one has
$u_\varepsilon(0) = \bar u(0) + \psit(0) - \psio(0)$.
\end{remark}
\begin{remark}\label{rem}
The rate ${1\over8}$ is not optimal. By adjusting the parameters~$\alpha$,
$\gamma$ and $\chi$ in an optimal way,
and by sharpening some of the arguments in our proof, one could achieve
a slightly better rate.
However, we do not believe that any rate obtained in this way would
reflect the true speed of convergence, so we keep
with the values (\ref{eexponents}) that yield simple fractions.
\end{remark}
\begin{remark}
From a~technical point of view, the general methodology followed in
this section and the subsequent sections is inspired from~\cite{Hai10},
where a~somewhat similar phenomenon was investigated.
Besides the structural differences in the equations considered here and
in \cite{Hai10},
the main technical difficulties that need to be overcome for the
present work
are the following:
\begin{longlist}[(1)]
\item[(1)] In \cite{Hai10}, it is possible to simply subtract the
stochastic convolution $\psi$ (or~$\widetilde\psi$) and work with the
equation for the remainder. Here, we instead subtract only the highest
Fourier modes of $\psi$. The reason for this choice
is that it entails that $\bar v^\gamma\to\bar u$ as $\varepsilon\to
0$. This allows us to linearize the nonlinearity around $\bar v^\gamma$
in order to exhibit the desired correction term. As a~consequence, our
a priori regularity estimates are much worse than those
in \cite{Hai10} and our convergence rates are worse. The main reason
why we need this complication is that our approximate derivative
$D_\varepsilon$ does not satisfy
the chain rule.
\item[(2)] All of our fixpoint arguments need to be performed in the
fractional Sobolev space $H^\alpha$, for some $\alpha> \frac12$.
This is in contrast to \cite{Hai10} where some of the arguments could
be performed first in $L^\infty$, and then lifted to $H^\alpha$ by a~standard bootstrapping argument. These bootstrapping arguments fail
here, since the nonlinearity of our approximating equation contains an
approximate derivative, which gives rise to correction terms which are
not easy to control.
\item[(3)] In one crucial step where a~Gaussian concentration
inequality is
employed in \cite{Hai10}, it was necessary that the
stochastic convolutions belong to~$H^\alpha$ for some $\alpha> \frac
12$. This is the case in \cite{Hai10} as a~consequence of the
extra regularizing effect caused by a~small fourth-order term present
in the linear part. This additional regularizing effect is not always present
in the current work. We therefore perform another truncation in Fourier
space, at very high frequencies. This is the
purpose of the exponent $\chi$.\vadjust{\goodbreak}
\end{longlist}
Note also that Proposition~\ref{propsecond-chaos} is the analogue of
Proposition 4.1 in \cite{Hai10}. One
difference is that we have a~much cleaner separation of the
probabilistic and the analytical aspects of this result.
\end{remark}

By a~standard Picard fixed point argument (see, e.g., \cite{Hai09}) it
can be shown that (\ref{equ-0}) admits a~unique mild
solution $\bar u$ defined on a~random time interval $[0,\tau_*]$.
Moreover, the spatial regularity of $\psio$ and $\bar u$
equals that of a~Brownian path, in the sense that $\psio(t)$ and $\bar
u(t)$ are continuous and belong to $H^\beta$ for
any $\beta< \frac12$ and any $t > 0$, but not to $H^{1/2}$. We
shall take advantage of the fact that the process $\vb$
is much more regular. In fact, $\vb(t)\in H^\beta$ almost surely for
any $\beta< \frac32$ and any $t > 0$, but one
does not expect it to belong to~$H^{3/2}$ in general. This follows
immediately from the mild formulation~(\ref{equ-0})
combined with a~standard bootstrapping argument.
It follows from these considerations that, for every fixed time
horizon~$T$, the stopping time
\[
\tau_*^K := T \wedge\inf\{ t\dvtx
\| \vb(t) \|_{\alpha}
\vee
\| \bar u(t) \|_{L^\infty}
\geq K \}
\]
converges in probability to $\tau_*\wedge T$ as $K \to\infty$.

It will be shown in Section~\ref{secprobtools} that a~number of
functionals of $\psio$ and $\psit$ scale in the following way:
\begin{eqnarray*}
\|\psitgx(t) \|_{L^\infty} &\lesssim& \varepsilon^{{\gamma}/2 -
\kappa} ,\qquad
\|\psitg(t) \|_{L^\infty} \lesssim\varepsilon^{{\gamma}/2 -
\kappa} ,\\
\|\psiog(t) \|_{L^\infty}
&\lesssim& \varepsilon^{{\gamma}/2 - \kappa} ,\qquad
\| \psitlx(t) \|_{L^\infty} \lesssim\varepsilon^{{\chi}/2 -
\kappa} ,\\
\|\psioug(t) \|_{\alpha}
&\lesssim&\varepsilon^{-\gamma(\alpha- {1/2}) - \kappa} ,\qquad
\| \psitgx(t)\|_{\alpha}
\lesssim\varepsilon^{-\chi(\alpha- 1/2) - \kappa} ,\\
\| \psitug(t) - \psioug(t)\|_{\alpha}
&\lesssim&\varepsilon^{2 - \gamma(\alpha+ 3/2) - \kappa} ,\qquad
\Theta_\varepsilon(\psitg(t))
\lesssim\varepsilon^{-1 - \kappa} ,\\
\Theta_\varepsilon(\psitlx(t))
&\lesssim&\varepsilon^{\chi- 2 - \kappa} ,\qquad
\|\Lambda- \Xi_\varepsilon(\psitgx(t))\|_{-\alpha} \lesssim
\varepsilon^{{1/2}-\kappa} ,
\end{eqnarray*}
where the quantities $\Theta_\varepsilon$ and $\Xi_\varepsilon$ are
defined by
\[
\Theta_\varepsilon(u) \stackrel{\mathrm{def}}{=}
\int_\R y^2 \| \hD_{\varepsilon y} u \|_{L^2}^2
|\mu|(dy) ,\qquad
\Xi_\varepsilon(u) \stackrel{\mathrm{def}}{=}\int_{\R} \frac
{\varepsilon y^2}{2}
\hD_{\varepsilon y} u \ot\hD_{\varepsilon y} u \mu(dy) .
\]
Note that all of these relations are of the form $\bolds{\Psi
}_i^\varepsilon(t) \lesssim\varepsilon^{\alpha_i-\kappa}$ for some
expression $\bolds{\Psi}_i^\varepsilon$
depending on $\varepsilon$ and some exponent $\alpha_i$.
In the proof, it will be convenient to impose this behavior by means
of a~hard constraint. For this purpose, we introduce the stopping time
$\tau^K$, which is defined for $K > 0$ by
%
%
\begin{equation} \label{eqtau-K}
\tau^K \stackrel{\mathrm{def}}{=}\tau_*^K \wedge\inf\{ t\dvtx\exists
i\dvtx
\bolds{\Psi}_i^\varepsilon(t) \geq\varepsilon^{\alpha_i -
\kappa} \} .
\end{equation}
From now on, we will write $C_K$ to denote a~constant which may depend
on~$K$ (and $T$) and is allowed to change
from line to line. Similarly, $\kappa$ will be a~positive universal
constant which is sufficiently small and whose value is allowed to
change from line to line. However, the final value of $\kappa$ is
independent of~$\varepsilon$, $K$ and $T$.

The remainder of this section is devoted to the proof of
Theorem~\ref{thmmain}.\vadjust{\goodbreak}
\begin{pf*}{Proof of Theorem~\ref{thmmain}}
Most of the work in the proof consists of bounding the difference
between $\vtg$ and $\vbg$ in $H^\alpha$. This\vspace*{1pt}
bound will be obtained in several steps, using the intermediate
processes $v_\varepsilon^{(i)}$, $i = 1, \ldots, 4$, defined by
\begin{subequation}\label{ev}
\begin{eqnarray}\qquad\qquad
\label{eqv-1}
v_\varepsilon^{(1)}(t) & = & S(t) v_0 + \psioug(t)
+ S * \bigl( \tF\bigl(v_\varepsilon^{(1)}\bigr) + \partial_x G\bigl
(v_\varepsilon^{(1)}\bigr)
\bigr)(t) ,\\
\label{eqv-2}
v_\varepsilon^{(2)}(t) & = & S(t) v_0 + \psioug(t)
+ S * \bigl( \tF\bigl(v_\varepsilon^{(2)}\bigr) + D_\varepsilon G\bigl
(v_\varepsilon^{(2)}\bigr)
\bigr)(t) ,\\
\label{eqv-3}
v_\varepsilon^{(3)}(t) & = & S_\varepsilon(t) v_0 + \psitug(t)
+ S_\varepsilon* \bigl( \tF\bigl(v_\varepsilon^{(3)}\bigr) +
D_\varepsilon
G\bigl(v_\varepsilon^{(3)}\bigr)
\bigr)(t) ,\\
\label{eqv-4}
v_\varepsilon^{(4)}(t) & = & S_\varepsilon(t) v_0
+ \psitug(t)
\nonumber\\[-8pt]\\[-8pt]
&&{} + S_\varepsilon* \bigl( F\bigl(v_\varepsilon^{(4)} + \psitgx\bigr)
+ \nabla G\bigl(v_\varepsilon^{(4)}
+ \psitgx\bigr)
D_\varepsilon\bigl(v_\varepsilon^{(4)} + \psitgx\bigr) \bigr)(t)
.\nonumber
\end{eqnarray}
\end{subequation}
At this stage, we stress that the main difficulty of the proof consists
of showing that $v_\varepsilon^{(3)}$ and $v_\varepsilon^{(4)}$ are
close (see Proposition~\ref{propest-34} below).\vspace*{1pt} Showing
the smallness
of the remaining differences $v_\varepsilon^{(j)}-v_\varepsilon
^{(j+1)}$ is relatively straightforward and follows by applying
standard SPDE techniques.
The main\vspace*{1pt} ingredient in this part of the proof\vspace*{1pt} is an
estimate which
compares the square of the approximate derivative of $\psitgx$ to the
correction term, in a~suitable Sobolev space of negative order. The
estimate is purely probabilistic and ultimately relies on the fact that
the quantity that we wish to control belongs to the second order Wiener
chaos. It can be found in Proposition~\ref{propsecond-chaos}, which
we consider to be the core of the paper.\vspace*{1pt}

Recall the definition of the stopping time $\tau^K$ given in (\ref
{eqtau-K}). With this definition
at hand, we set $\tau^K_0 \stackrel{\mathrm{def}}{=}\tau^K$ as well as
$v_\varepsilon
^{(0)} \stackrel{\mathrm{def}}{=}\vbg$ and $v_\varepsilon^{(5)}
\stackrel
{\mathrm{def}}{=}\vtg,$ and we
define recursively a~sequence of stopping times
$\tau^K_j$ with $j=1,\ldots,5$ by
%
%
\begin{equation}\label{edeftauKj}
\tau^K_j = \tau^K_{j-1} \wedge\inf\bigl\{t\dvtx\bigl\|v_\varepsilon
^{(j)}(t)-v_\varepsilon^{(j-1)}(t)\bigr\|_\alpha\ge K\bigr\}.
\end{equation}
With this notation at hand, Propositions~\ref{propest-01}--\ref
{propest-4eps} state that,
for all fixed values $K, \kappa>0$ and every $j = 1,\ldots,5$, one has
%
%
\begin{equation}\label{eboundvvj}
\lim_{\varepsilon\to0} {\mathbb P}\Bigl( \sup_{t \leq\tau_j^K} \bigl\|
v_\varepsilon^{(j)}(t) - v_\varepsilon^{(j-1)}(t) \bigr\|_{\alpha} >
\varepsilon^{1/8 - \kappa} \Bigr) = 0 .
\end{equation}
Combining all of these bounds, we conclude immediately that, for every
fixed time horizon $T>0$ and every
choice of values $K$ and $\kappa$, we have
\[
\lim_{\varepsilon\to0} {\mathbb P}\Bigl( \sup_{t \leq\tau_5^K} \|
\vtg(t) - \vbg(t) \|_{\alpha} > \varepsilon^{1/8 - \kappa}
\Bigr) = 0 .
\]
This is formally very close to (\ref{eboundvv}), except that we still
have the values $T, K>0$ appearing in our statement
and consider the solutions only up to the stopping time~$\tau_5^K$.

Since $\tau_* \wedge T \to\tau_*$ as $T\to\infty$ and since we
already argued that $\tau_*^K \to\tau_* \wedge T$
as $K\to\infty$, the bound (\ref{eboundvv}) follows if we are
able\vadjust{\goodbreak}
to show that, for every fixed choice of $K$ and $T$,
%
%
\begin{equation}\label{etauequ}
\lim_{\varepsilon\to0} {\mathbb P}( \tau_5^K \neq\tau_*^K ) = 0 .
\end{equation}
Since\vspace*{1pt} the statement of our theorem is stronger, the smaller the value
of $\kappa$, we can assume without
loss of generality that $\kappa< {1\over8}$. In this case, $\lim
_{\varepsilon\to0}\varepsilon^{{1/8}-\kappa} = 0$, so that
(\ref{eboundvvj}) and (\ref{edeftauKj}) together imply that
\[
\lim_{\varepsilon\to0} {\mathbb P}( \tau_j^K \neq\tau_{j-1}^K ) =
0
\]
for $j=1,\ldots,5$, from which we conclude that $\lim_{\varepsilon
\to0} {\mathbb P}( \tau_5^K \neq\tau^K ) = 0$.

In order to finish the proof of (\ref{eboundvv}), it now suffices to
show that\break
$\lim_{\varepsilon\to0} {\mathbb P}( \tau^K \neq\tau_*^K ) = 0$.
Fix an arbitrary $T>0$ and $\kappa> 0$.
It then follows from Propositions~\ref{proppsi-est}, \ref
{proppsi-comparison} and~\ref{propDpsi-est} that
for each of the terms $\bolds{\Psi}_j^\varepsilon$ appearing in (\ref{eqtau-K}),
there exists a~constant $C_j>0$ such that
\[
\E\sup_{t \in[0,T]} \bolds{\Psi}_j^\varepsilon(t) \leq C_j
\varepsilon^{\alpha_j - {\kappa/2}} ,
\]
uniformly for all $\varepsilon\le1$. It then follows from Chebychev's
inequality that
\[
{\mathbb P}( \tau^K \neq\tau_*^K ) \le\sum_j {\mathbb P}
\Bigl(\sup_{t \in[0,T]} \bolds{\Psi}_j^\varepsilon(t) \ge
\varepsilon^{\alpha_j-\kappa} \Bigr)
\le\sum_j C_j \varepsilon^{\kappa/ 2} ,
\]
from which the claim follows.

Since (\ref{ebounduu}) follows from (\ref{eboundvv}) and (\ref
{eboundpsipsi}), the proof of the theorem is complete if we show
that (\ref{eboundpsipsi}) holds.
Since it follows from Proposition~\ref{proppsi-est} and Chebychev's
inequality that
\[
\lim_{\varepsilon\to0} {\mathbb P}\Bigl( \sup_{t \leq T} \| \psitg
(t) - \psiog(t) \|_{L^\infty} > \varepsilon^{1/2 - \kappa}
\Bigr) = 0
\]
for every $T>0$, this claim follows at once.
\end{pf*}

\section{Analytic tools}
\label{secingredients}

\subsection{Products and compositions of functions in Sobolev spaces}

In this subsection, we collect some well-known bounds for products and
compositions of functions in Sobolev spaces.
As is usual in the analysis literature, we use the notation $\Phi
\lesssim\Psi$ as a~shorthand for ``there exists a~constant $C$ such
that $\Phi\le C \Psi$.'' These estimates will be useful in order
control the various terms that arise in the Taylor expansion of the
nonlinearity that will be performed in Section~\ref{secapprox} below.
\begin{lemma} \label{lemproduct-sobolev}
Let $r,s,t \geq0$ be such that $r \wedge s > t$ and $r + s > \frac12 +t$.
\begin{longlist}[(1)]
\item[(1)]
For $f \in H^r$ and $g \in H^s$, we have $fg \in H^t$ and
%
%
\begin{equation} \label{eqSobolev-product}
\| fg \|_t \lesssim\| f \|_r \| g \|_s .
\end{equation}
\item[(2)]
For $f \in H^r$ and $g \in H^{-t}$, we have $fg \in H^{-s}$ and
%
%
\begin{equation} \label{eqSobolev-product-neg}
\| fg \|_{-s}
\lesssim\| f \|_{r} \| g \|_{-t} .
\end{equation}
\end{longlist}
\end{lemma}
\begin{pf}
This result is very well known. A proof of (\ref{eqSobolev-product})
can be found, for example, in \cite{Hai09}, Theorem 6.18, and (\ref
{eqSobolev-product-neg}) follows by duality.
\end{pf}
\begin{lemma}\label{lemsobolev-composition}
Let $s \in(\frac12,1)$. There exists $C > 0$ such that for any $u \in
H^s$ and any $G \in C^1(\R^n;\R^n)$ satisfying
\[
\| G_u \|_{C^1} := \sup\{ | G(x) |+ | \nabla G(x)|
\dvtx|x| \leq\| u \|_{L^\infty} \} < \infty,
\]
we have
\[
\|G \circ u\|_s \leq C \| G_u \|_{C^1} ( 1 + \|u\|_s) .
\]
\end{lemma}
\begin{pf}
Let $\tau_h$ be the shift operator defined by $\tau_h u(x) :=
u(x-h)$. It is well known (see, e.g., \cite{DVS93} or, for functions
defined on $\R^n$, \cite{AF03}, Theorem 7.47) that the expression
%
%
\begin{equation}\label{eequivHs}
\| u \|_{L^2} +
\biggl( \int_0^1 \Bigl[ t^{-s} \sup_{|h| < t}
\| u - \tau_h u \|_{L^2} \Bigr]^2 \,\frac{dt}{t} \biggr)^{1/2}
\end{equation}
defines an equivalent norm on $H^s$. The result then follows by
inserting the estimates
\begin{eqnarray*}
\| G \circ u \|_{L^2} & \leq & \| G \circ u \|_{L^\infty}
\leq C \| G_u \|_{C^1},\\
\| G \circ u - \tau_h (G \circ u) \|_{L^2}
& \leq & C\| G_u \|_{C^1} \| u - \tau_h u \|_{L^2}
\end{eqnarray*}
into (\ref{eequivHs}).
\end{pf}

\subsection{Semigroup bounds}

We will frequently use the fact that for $\alpha\geq\beta$ and $T >
0$, there exists a~constant $C > 0$ such that
%
%
\begin{equation} \label{eqS-bound}
\| S(t) u \|_{\alpha} \leq C t^{-(\alpha- \beta)/2} \| u \|
_{\beta}
\end{equation}
for any $\varepsilon\in(0,1]$, $t \in[0,T]$ and $u \in H^\beta$.
This is a~straightforward consequence of
standard analytic semigroup theory \cite{Lunardi,Hai09}. Since the
generator of $S$ is selfadjoint in all of the $H^s$, it is
also straightforward to prove (\ref{eqS-bound}) by hand. As a~consequence, we have:
\begin{lemma} \label{lemsemigroup-reg}
Let $\alpha, \beta\in\R$ be such that $0\leq\alpha- \beta< 2$
and let $T > 0$. There exists $C > 0$ such that for all $t \in[0,T]$
and $u \in C([0,t]; H^\beta)$ we have
%
%
\begin{equation} \label{eqsemigroup-reg}
\biggl\| \int_0^t S(t-s) u(s) \,ds \biggr\|_{\alpha}
\leq Ct^{1 - (\alpha- \beta)/2}
\sup_{s \in[0,t]} \| u(s) \|_{\beta} .
\end{equation}
\end{lemma}
\begin{pf}
It suffices to integrate the bound (\ref{eqS-bound}).\vadjust{\goodbreak}
\end{pf}

The following bounds measure how well $S_\varepsilon$ approximates $S$
in these interpolation spaces.
The general philosophy is that every power of $\varepsilon$ has to be
paid with one spatial derivative worth
of regularity. This type of power-counting is a~direct consequence of
the fact that the function $f$ that measures
how much $\Delta_\varepsilon$ differs from $\partial_x^2$, is
evaluated at
$\varepsilon|k|$ in the definition of $\Delta_\varepsilon$.
The precise bounds are the following:
\begin{lemma}\label{lemS-comparison}
Let $\kappa\in[0,2]$. For $T > 0$ there exists $C > 0$ such that for
any $t \in[0,T]$, $\varepsilon\in(0,1]$, and $u \in H^\beta$, we have
%
%
\begin{eqnarray} \label{eqS-comparison-1}\qquad
\| S_\varepsilon(t) u - S(t) u \|_{\alpha}
& \leq & C \varepsilon^{\kappa} t^{-(\alpha- \beta+ \kappa)/2}
\| u \|_{\beta}  \qquad(\beta\leq\alpha
+ 2\kappa) , \\
\label{eqS-comparison-2}
\| S_\varepsilon(t) u \|_{\alpha}
& \leq & C t^{-(\alpha- \beta)/2} \| u \|_{\beta} \qquad(\beta\leq
\alpha) .
\end{eqnarray}
\end{lemma}
\begin{pf}
We set\vspace*{1pt} $\bar f \stackrel{\mathrm{def}}{=}f-1$ and assume $\nu= 1$ for notational
simplicity, since the case
$\nu\neq1$ is virtually identical. The assumptions on $f$ imply that
$
|\bar f(\varepsilon n)| \leq c \varepsilon^2 n^2
$
whenever $n < \delta/ \varepsilon$ and $\delta$ is some sufficiently
small constant. Using the mean value theorem and the fact that
we can assume $\delta<1$ without loss of generality, we obtain for $n
< \delta/ \varepsilon$ and $\kappa\in[0,2]$,
\begin{eqnarray*}
|{\exp}(-tn^2 \bar f(\varepsilon n) ) - 1|
& \leq &(2 \wedge c t \varepsilon^2 n^4 )e^{c t \varepsilon^2 n^4}
\leq C t^{{\kappa}/{2}} \varepsilon^{\kappa} n^{2\kappa}e^{c t
\delta^2 n^2}
\\
&\leq& C \varepsilon^{\kappa} t^{{\kappa}/{2}} n^{2\kappa
}e^{c \delta^2 t (1+n^2)} .
\end{eqnarray*}
Inserting this bound into the identity
\[
\bigl(S_\varepsilon(t) u - S(t)\bigr) e_n
= \bigl(e^{-tn^2 \bar f(\varepsilon n)} - 1\bigr) e^{-t(1+n^2)} e_n ,
\]
it then follows from (\ref{eqS-bound}) that
\begin{eqnarray}
\label{eqlow-freq}
\bigl\| \Pi_{\delta/\varepsilon}\bigl(S_\varepsilon(t) - S(t)\bigr) u \bigr\|_{\alpha}
& \leq & C \varepsilon^{\kappa} t^{{\kappa}/{2}} \bigl\| S\bigl((1-\delta
^2 c)t\bigr) u \bigr\|_{\alpha+ 2\kappa}
\nonumber\\[-8pt]\\[-8pt]
& \leq & C \varepsilon^{\kappa} t^{-(\alpha- \beta+ \kappa)/2}
\| u \|_{\beta} ,
\nonumber
\end{eqnarray}
provided that we choose $\delta$ sufficiently small so that $\delta^2
c \le{1\over2}$, say.

On the other hand, note that
\[
(I - \Pi_{\delta/\varepsilon})\bigl(S_\varepsilon(t) u - S(t)\bigr) e_n
= \one_{\{ | n | > \delta/ \varepsilon\}}\bigl(e^{-tn^2 \bar
f(\varepsilon n)} - 1\bigr) e^{-t(1+n^2)} e_n .
\]
Recall that $\bar f(\varepsilon n) \geq q-1$ for all $n$, and that $q
\in(0,1]$. Then we can find a~constant $C$ such that
\[
|{\exp}(-tn^2 \bar f(\varepsilon n) ) - 1| e^{-t(1+n^2)}
\leq C e^{-q(1+n^2) t} .
\]
Moreover, for any $\kappa> 0$ we have $\one_{\{ | n | > \delta/
\varepsilon\}} \leq|\varepsilon n/\delta|^{\kappa}$.
It thus follows, using~(\ref{eqS-bound}) again, that
\[
\bigl\| ( I - \Pi_{\delta/\varepsilon})\bigl(S_\varepsilon(t) - S(t)\bigr) u \bigr\|
_{\alpha}
\leq C \varepsilon^{\kappa} \| S(qt) u \|_{\alpha+ \kappa}
\leq C \varepsilon^{\kappa}
t^{-(\alpha- \beta+ \kappa)/2} \| u \|_{\beta} .
\]
The bound (\ref{eqS-comparison-1}) now follows by combining this
inequality with (\ref{eqlow-freq}).
Inequality (\ref{eqS-comparison-2}) follows by combining the
special case $\kappa= 0$ with (\ref{eqS-bound}).
\end{pf}

\subsection{Estimates for the gradient term}

In this section, we similarly show how well the operator $D_\varepsilon
$ approximates $\partial_x$. Again, the guiding
principle is that every power of $\varepsilon$ ``costs'' the loss of
one derivative. However, we are also going to use
the fact that $D_\varepsilon$ is a~bounded operator. In this case, we
can gain up to one spatial derivative with respect to
the operator $\partial_x$, but we have to ``pay'' with the same number of
inverse powers of $\varepsilon$. The rigorous statement
for the latter fact is the following lemma.
\begin{lemma}\label{lemD-eps-grad}
Let $\beta\in\R$ and $\alpha\in[0,1]$. There exists $C > 0$ such
that for all $\varepsilon\in(0,1]$ and $u \in H^\beta$ the estimate
\[
\| D_\varepsilon u\|_{\beta-\alpha} \leq C \varepsilon^{\alpha-1}
\| u\|_\beta
\]
holds.
\end{lemma}
\begin{pf}
Using the assumption that $M := |\mu|(\R) < \infty$, together with
Jensen's inequality and Fubini's theorem, we obtain
\begin{eqnarray*}
\| D_\varepsilon u \|_{L^2}^2
& \leq &\frac{1}{\varepsilon^2}
\int\biggl( \int_\R|u(x+ \varepsilon y)| |\mu|(dy) \biggr)^2 \, dx
\\& \leq &\frac{M}{\varepsilon^2}
\iint_\R|u(x+ \varepsilon y)|^2 |\mu|(dy) \, dx = \frac
{M^2}{\varepsilon^2} \| u \|_{L^2}^2 .
\end{eqnarray*}

On the other hand, assuming for the moment that $u$ is smooth, we use
the assumption that $\mu(\R) = 0$, and apply Jensen's inequality and
Minkowski's integral inequality to obtain
\begin{eqnarray*}
\| D_\varepsilon u \|_{L^2}^2
& = & \frac{1}{\varepsilon^2}
\int\biggl( \int_\R u(x + \varepsilon y) \mu(dy) \biggr)^2 \, dx
\\
& = &\frac{1}{\varepsilon^2}
\int\biggl( \int_\R\int_0^{\varepsilon y} u'(x + z) \,dz \,\mu
(dy) \biggr)^2 \, dx
\\
& \leq &\frac{M}{\varepsilon^2}
\iint_\R\biggl( \int_0^{\varepsilon y} |u'(x + z)| \,dz
\biggr)^2 |\mu|(dy) \, dx
\\
& \leq &\frac{M}{\varepsilon^2}
\int_\R\biggl( \int_0^{\varepsilon y} \biggl( \int|u'(x + z)|^2
\, dx \biggr)^{1/2} \,dz \biggr)^2 |\mu|(dy)
\\
& = &M \|u' \|_{L^2}^2
\int_\R y^2 |\mu|(dy) \leq C \|u\|_1^2 .
\end{eqnarray*}
Using complex interpolation, it follows that $\| D_\varepsilon u \|
_{L^2} \leq C\varepsilon^{\alpha- 1} \| u \|_\alpha$ for every
$\alpha\in[0,1]$.
The desired result then follows from the fact that $D_\varepsilon$
commutes with every Fourier multiplier.
\end{pf}

The announced approximation result on the other hand is the following lemma.
\begin{lemma}\label{lemD-eps-bound}
Let $\beta\in\R$ and $\alpha\in[0,1]$. There exists $C > 0$ such
that for all $\varepsilon\in(0,1]$ and $u \in H^\beta$ the estimate
\[
\| D_\varepsilon u - \partial_x u \|_{\beta-1-\alpha} \leq C
\varepsilon^\alpha\| u \|_\beta
\]
holds.
\end{lemma}
\begin{pf}
In view of (\ref{eqmoments}) we have, assuming for the moment that
$u$ is smooth,
\[
(D_\varepsilon- \partial_x) u(x)
= \frac{1}{\varepsilon} \int_{\R} \int_0^{\varepsilon y} \int_0^w
u''(x + z) \,dz \,dw \,\mu(dy) .
\]
Integrating against a~test function $\varphi$ and applying Fubini's
theorem, we arrive at
\begin{eqnarray*}
\biggl| \int\varphi(x)(D_\varepsilon- \partial_x) u(x) \, dx \biggr|
& \leq &\frac{C}{\varepsilon} \int_{\R} \int_0^{\varepsilon y} \int_0^w
\|\varphi\|_{2-\beta} \|u\|_{\beta} \,dz \,dw\, |\mu|(dy)
\\&\leq &C\varepsilon\|\varphi\|_{2-\beta} \|u\|_{\beta}\int_{\R}
|y|^2 |\mu|(dy) ,
\end{eqnarray*}
which implies that
\[
\| (D_\varepsilon- \partial_x) u \|_{\beta-2} \leq C \varepsilon\|
u \|_{\beta} .
\]
On the other hand, Lemma~\ref{lemD-eps-grad} implies that
\[
\| ( D_\varepsilon- \partial_x) u \|_{\beta-1} \leq C \|u\|_{\beta
} ,
\]
and the result then follows as before by interpolating between these estimates.
\end{pf}

As an immediate corollary of these bounds, we obtain the following
useful fact.
\begin{corollary}\label{corD-psi-eps-bound}
Let $\beta\in[0,1)$. There exists $C > 0$ such that for $\varepsilon
\in(0,1]$, $u \in H^{\beta}$, and $G \in C^1(\R^n)$ we have
\[
\| D_\varepsilon G(u) - \partial_x G(u) \|_{-1}
\leq C \varepsilon^\beta\|G_u\|_{C^1}
( 1 + \| u \|_\beta) ,
\]
where $\|G_u \|_{C^1}$ is defined as in Lemma~\ref{lemsobolev-composition}.
\end{corollary}
\begin{pf}
Using Lemmas~\ref{lemD-eps-bound} and~\ref{lemsobolev-composition},
we obtain
\[
\| D_\varepsilon G(u) - \partial_x G(u) \|_{-1}
\leq C\varepsilon^{\beta}\| G(u) \|_\beta
\leq C\varepsilon^\beta\|G_u\|_{C^1} ( 1 + \| u \|_\beta) ,
\]
which is the stated claim.
\end{pf}

\section{Probabilistic tools}
\label{secprobtools}

In this section, we prove some sharp estimates for certain expressions
involving stochastic convolutions. Our main tool is the following\vadjust{\goodbreak}
version of Kolmogorov's continuity criterion, which follows immediately
from the one given, for example, in \cite{RevYor}.
The reason why we state condition (\ref{eqequiv-mom}) in this form,
is that it is automatically satisfied (by hypercontractivity) for
random fields taking values in a~Wiener chaos of fixed (finite) order.
\begin{lemma}\label{lemKolmogorov}
Let $(\varphi(t))_{t \in[0,1]^n}$ be a~Banach space-valued random
field having the property that for any $q \in(2,\infty)$ there exists
a constant $K_q > 0$ such that
\begin{eqnarray}
\label{eqequiv-mom}
(\E\| \varphi(t) \|^q )^{1/q} & \leq & K_q
(\E\| \varphi(t) \|^2 )^{1/2} ,\nonumber\\[-8pt]\\[-8pt]
\bigl(\E\| \varphi(s) - \varphi(t) \|^q \bigr)^{1/q} & \leq & K_q
\bigl(\E\| \varphi(s) - \varphi(t) \|^2 \bigr)^{1/2}
\nonumber
\end{eqnarray}
for all $s, t \in[0,1]^n$. Furthermore, suppose that the estimate
\[
\E\| \varphi(s) - \varphi(t) \|^2 \leq K_0 | s - t |^\delta
\]
holds for some $K_0, \delta> 0$ and all $s,t \in[0,1]^n$. Then, for
every $p > 0$ there exists $C > 0$ such that
\[
\E\sup_{t \in[0,1]^n} \| \varphi(t) \|^p
\leq C \bigl( K_0 + \E\|\varphi(0)\|^2 \bigr)^{{p/2}} .
\]
\end{lemma}

Throughout\vspace*{1pt} this subsection, we shall use $\thk$ and $\thtk$ for the
Fourier coefficients of
$\psio$ and $\psit$, so that
\[
\psio(t) = \sum_{k \in\Z}\thk(t) e_k ,\qquad
\psit(t) = \sum_{k \in\Z}\thtk(t) e_k .
\]
With this notation at hand, we first state the following approximation
bound, which shows that we can
again trade powers of $k$ for powers of~$\varepsilon$, provided that
we look at the difference squared:
\begin{lemma}\label{lemclaim}
For $t \geq0$, $k \in\Z$ and $\varepsilon\in(0,1]$, we have
%
%
\begin{equation} \label{eqclaim}
\E| \thtk(t) - \thk(t) |^2
\leq C (k^{-2} \wedge\varepsilon^4 k^2 ) .
\end{equation}
\end{lemma}
\begin{pf}
We write again $\bar f = f-1$ and assume $\nu= 1$ for simplicity. The
It\^{o} isometry then implies that
\begin{eqnarray}
\label{eqIto}
\E|\thtk(t) - \thk(t)|^2
& = & C\int_0^\infty e^{-2t(1+k^2)} \bigl(
1 - h(\varepsilon|k|)e^{-t k^2 \bar f(\varepsilon|k|)} \bigr)^2 \,dt
\nonumber\\
& \leq & C \int_0^\infty e^{-2t(1+k^2)} \bigl(
1 - e^{-t k^2 \bar f(\varepsilon|k|)} \bigr)^2 \,dt
\nonumber\\[-8pt]\\[-8pt]
&&{} + C \int_0^\infty e^{-2t(1+k^2)}
e^{-2 t k^2 \bar f(\varepsilon|k|)}
\bigl( 1 - h(\varepsilon|k|)\bigr)^2 \,dt
\nonumber\\
& \stackrel{\mathrm{def}}{=} &I_1 + I_2 .\nonumber
\end{eqnarray}
Let $\delta> 0$ be a~(small) constant to be determined later
and consider first the term $I_1$ with $|\varepsilon k| \leq\delta$.
Since $f$ is twice differentiable near the origin, we can find $\delta
$ small enough so that
$|\bar f(|\varepsilon k|)| \leq c \varepsilon^2 k^2$ for some $c>0$.
Therefore, for $t \geq0$,
%
%
\begin{equation} \label{eqexp-est}
\bigl|1 - e^{-t k^2 \bar f(\varepsilon|k|)}\bigr|
\leq c t \varepsilon^2 k^4 e^{c t \varepsilon^2 k^4} \le c t
\varepsilon^2 k^4 e^{c \delta^2 t k^2} ,
\end{equation}
so that
\[
|I_1| \le C \varepsilon^4 k^8 \int_0^{\infty} t^2 e^{-2t(1+k^2)+ 2c
\delta^2 k^2 t} \,dt .
\]
If we ensure that $\delta$ is small enough so that $2c \delta^2 \le
1$, we obtain
\[
|I_1| \leq C \varepsilon^4 k^8
\int_0^{\infty} t^2 e^{-k^2 t} \,dt
\leq C \varepsilon^4 k^2 \le C (k^{-2} \wedge\varepsilon^4
k^2) ,
\]
where the last inequality follows from the fact that $|\varepsilon k|
\leq\delta$ by assumption.

To treat the case $|\varepsilon k| > \delta$,
we use the fact that by assumption there exists $q \in(0,1]$ such that
$f \ge q$, so that
\begin{eqnarray} \label{eqclaim-2}
|I_1| &\le& \int_{0}^\infty e^{-2t k^2}
\bigl( 1 - e^{-t k^2 (q - 1)} \bigr)^2 \,dt
\le C \int_0^\infty e^{-2t q k^2} \,dt \nonumber\\[-8pt]\\[-8pt]
&\le& C k^{-2} \le C (k^{-2} \wedge\varepsilon^4
k^2) .\nonumber
\end{eqnarray}

The bound on $I_2$ works in pretty much the same way, using the fact
that the assumptions
on $h$ imply that
\[
\bigl|1-h(\varepsilon|k|)\bigr| \le C (1 \wedge\varepsilon^2 k^2) .
\]
Using again the fact that $f \ge q$, we then obtain
\[
I_2 \le C\int_0^\infty e^{-2t q k^2}(1 \wedge\varepsilon^4
k^4) \,dt \le C (k^{-2} \wedge\varepsilon^4 k^2)
\]
as required.
\end{pf}

We continue with a~sequence of propositions, in which the estimates
obtained in the previous lemma are used to establish various bounds for
stochastic convolutions.
\begin{proposition}\label{proppsi-est}
Let $0 < \gamma< \chi$. For $\kappa> 0$ and $\varepsilon\in(0,1]$
we have
\begin{eqnarray*}
\E\sup_{t \in[0,T]} \| \psiog(t) \|_{L^\infty}
&\leq& C \varepsilon^{{\gamma}/2 - \kappa} ,\qquad
\E\sup_{t \in[0,T]} \| \psitg(t) \|_{L^\infty}
\leq C \varepsilon^{{\gamma}/2 - \kappa} , \\
\E\sup_{t \in[0,T]} \| \psitgx(t) \|_{L^\infty}
&\leq& C \varepsilon^{{\gamma}/2 - \kappa} ,\qquad
\E\sup_{t \in[0,T]}
\| \psitg(t) - \psiog(t)\|_{L^\infty}
\leq C \varepsilon^{1/2 - \kappa} .
\end{eqnarray*}
\end{proposition}
\begin{pf}
We start with the proof of the second estimate.
Observe that $\thtk$ is a~complex one-dimensional stationary
Ornstein--Uhlenbeck process with variance $h^2/(2(1 + \nu k^2 f))$ and
characteristic time $1 + \nu k^2 f$. This implies that
%
%
\begin{equation} \label{eqpsi-1}
\E| \thtk(t) |^2
= \frac{h^2(\varepsilon|k|)}{2(1 + \nu k^2 f(\varepsilon|k|))}
\leq C (1 \wedge k^{-2})
\end{equation}
and
%
%
\begin{equation} \label{eqpsi-2}
\E| \thtk(t) - \thtk(s) |^2
\leq C h^2(\varepsilon|k|)| t - s |
\leq C| t - s | .
\end{equation}
These bounds imply that, on the one hand,
\[
\E| \thtk(t) e_k(x) - \thtk(s) e_k(y) |^2
\leq C \E|\thtk(t)|^2 + C \E|\thtk(s)|^2
\leq C (1 \wedge k^{-2}) ,
\]
while on the other hand, one has
\begin{eqnarray*}
&& \E| \thtk(t) e_k(x) - \thtk(s) e_k(y) |^2
\\
&&\qquad \leq C \E| \thtk(t) - \thtk(s) |^2
+ C k^2 |x - y|^2 \E| \thtk(s)|^2
\\
&&\qquad \leq C | t - s | + C | x - y |^2 .
\end{eqnarray*}
Combining these inequalities we find that, for every $\kappa\in[0,2]$,
\[
\E| \thtk(t) e_k(x) - \thtk(s) e_k(y) |^2
\leq C (1 \wedge k^{-2})^{1 - \kappa/2} (| t - s | + | x - y
|^2)^{\kappa/2} .
\]
Since the $\thtk$'s are independent except for the reality condition
$\thtmk= \hspace*{1pt}\overline{\hspace*{-1pt}\widetilde{\theta}}_k$, we infer that
\begin{eqnarray*}
\E| \psitg(t,x) - \psitg(s,y) |^2
&\leq& C \sum_{|k| > \varepsilon^{-\gamma}} \E| \thtk(t) e_k(x) -
\thtk(s) e_k(y) |^2
\\
&\leq& C (| t - s | + | x - y |^2)^{\kappa/2}
\sum_{|k| > \varepsilon^{-\gamma}} (1 \wedge k^{-2})^{1 -
\kappa/2}
\\
&\leq& C \varepsilon^{(1 - \kappa)\gamma}
(| t - s | + | x - y |^2)^{\kappa/2} .
\end{eqnarray*}
Arguing similarly, we obtain
\[
\E|\psitg(0,0)|^2
\leq C \sum_{|k| > \varepsilon^{-\gamma}} \E|\thtk(0)|^2
\leq C \sum_{|k| > \varepsilon^{-\gamma}} (1 \wedge k^{-2})
\leq C \varepsilon^{\gamma} .
\]
The result now follows by combining these two bounds with Lemma \ref
{lemKolmogorov}.

The proof of the first and third estimates being very similar, we do
not reproduce them here.
In order to prove the last estimate, we use Lemma~\ref{lemclaim} to obtain
\[
\E| \thtk(t) - \thk(t)|^2 \leq C (k^{-2})^{3/4+\kappa/4}
(\varepsilon^4 k^2)^{1/4- \kappa/4}
\leq C \varepsilon^{1-\kappa}k^{-1-\kappa} .
\]
This bound then replaces (\ref{eqpsi-1}), and the rest of the proof
is again analogous
to the proof of the second estimate.
\end{pf}
\begin{proposition}\label{proppsi-comparison}
Let $\zeta> 0$. For $\kappa> 0$ and $\varepsilon\in(0,1]$, we have
%
%
\begin{eqnarray}
\label{eqpsi-comparison3}\qquad
\E\sup_{t \in[0,T]}
\| \psiouz(t) \|_{\alpha}
& \leq & C \varepsilon^{-\zeta(\alpha- 1/2) - \kappa}
\qquad\biggl(\alpha> \frac12\biggr) ,\\
\label{eqpsi-comparison1}
\E\sup_{t \in[0,T]}
\| \psituz(t) - \psiouz(t) \|_{\alpha}
& \leq & C \varepsilon^{2 - \zeta(\alpha+ 3/2) - \kappa}
\qquad\biggl(\alpha> - \frac32\biggr) .
\end{eqnarray}
\end{proposition}
\begin{pf}
In view of the estimates
%
%
\begin{equation} \label{eqtheta-bounds}
\E| \thk(t) |^2 \leq C k^{-2} ,\qquad
\E| \thk(t) - \thk(s) |^2 \leq C | t - s | ,
\end{equation}
we obtain
\begin{eqnarray*}
\E\| \psiouz(t) - \psiouz(s) \|_\alpha^2
& \leq& C |t-s|^\kappa
\sum_{|k| \leq\varepsilon^{-\zeta}} (1+ k^2)^{\alpha- 1 + \kappa}
\\& \leq& C |t-s|^\kappa
\varepsilon^{-2\zeta(\alpha- 1/2 +\kappa)}
\end{eqnarray*}
and
\[
\E\| \psiouz(0) \|_\alpha^2
\leq C \varepsilon^{-2\zeta(\alpha- 1/2)} .
\]
Inequality (\ref{eqpsi-comparison3}) thus follows from Lemma \ref
{lemKolmogorov}.

In order to prove (\ref{eqpsi-comparison1}), we argue similarly, but
the estimates are slightly more involved. Write $\delta_k := \thtk-
\thk$ so that $\psituz- \psiouz= \sum_{|k| \leq\varepsilon
^{-\zeta}} \delta_k e_k$. Using~(\ref{eqpsi-2}) and (\ref
{eqtheta-bounds}), we have for $s,t \geq0$,
\[
\E| \delta_k(t) - \delta_k(s) |^2
\leq C| t - s | .
\]
Combining this bound with Lemma~\ref{lemclaim}, we infer that for
$\kappa\in[0,\frac12)$,
\[
\E| \delta_k(t)
- \delta_k(s) |^2
\leq C (k^{-2})^\kappa
(\varepsilon^4 k^2)^{1-2\kappa} | t-s |^{\kappa}
= C \varepsilon^{4 - 8 \kappa} k ^{2 - 6 \kappa} | t-s |^{\kappa} .
\]
For $\kappa\in(0,\frac13 \alpha+ \frac12)$, we thus obtain
\begin{eqnarray*}
\E\| ( \psituz- \psiouz)(t)
- (\psituz- \psiouz)(s) \|_\alpha^2
& \leq & C |t-s|^\kappa\varepsilon^{4 - 8\kappa}
\sum_{|k| \leq\varepsilon^{-\zeta}} (1+ k^2)^{\alpha+ 1 - 3\kappa}
\\
& \leq & C |t-s|^\kappa
\varepsilon^{4 - \zeta(2\alpha+ 3) - 8\kappa}
\end{eqnarray*}
and similarly
\[
\E\sup_{t \in[0,T]}
\| \psituz(t) - \psiouz(t)\|_{\alpha}^2
\leq C \varepsilon^{4 - \zeta(2\alpha+ 3) - 8\kappa} .
\]
The desired estimate (\ref{eqpsi-comparison1}) now follows from Lemma
\ref{lemKolmogorov}.
\end{pf}
\begin{proposition}\label{propDpsi-est}
Let $\zeta> 0$. For every $\kappa> 0$ there exists $C>0$ such that
\[
\E\sup_{t\in[0,T]} \Theta(\psitz(t))
\leq C \varepsilon^{-1 + (\zeta-1)^+ - \kappa}
\]
for all $\varepsilon\in(0,1]$, where we wrote $(\zeta- 1)^+
\stackrel{\mathit{def}}{=}0
\vee(\zeta- 1)$.\vadjust{\goodbreak}
\end{proposition}
\begin{pf}
As in the proof of Propositions~\ref{proppsi-est} and \ref
{proppsi-comparison}, we shall apply Kolmogorov's continuity criterion
from Lemma~\ref{lemKolmogorov}, this time for $L^2$-valued random fields.
It follows from (\ref{eqpsi-1}) that
\begin{eqnarray*}
\E\bigl\| \hD_{\varepsilon y} \bigl(\psitz(t)
- \psitz(s) \bigr) \bigr\|_{L^2}^2
& = & \sum_{|k| > \varepsilon^{-\zeta}}
\biggl|\frac{e^{ik\varepsilon y} -1 }{\varepsilon y}\biggr|^2
\E| \thtk(t) - \thtk(s) |^2
\\
& \leq & C \sum_{k > \varepsilon^{-\zeta}}
\frac{1 - \cos(k\varepsilon y) }{|\varepsilon k y|^2} .
\end{eqnarray*}
Note that, up to a~factor $\varepsilon|y|$, this sum can be
interpreted as a~Riemann sum for the function $H(t) \stackrel
{\mathrm{def}}{=}t^{-2}(1 -
\cos(t))$. In fact, since $H(t) \leq2(1 \wedge t^{-2})$,
\begin{eqnarray}\label{eqRiemann-sum}\qquad
\varepsilon|y| \sum_{k > \varepsilon^{-\zeta}} \frac{1 - \cos
(k\varepsilon y)}{|k \varepsilon y|^2}
&=&
\sum_{k > \varepsilon^{-\zeta}} \varepsilon|y| H(k\varepsilon y)
\leq2
\int_{\varepsilon^{1-\zeta}}^\infty(1 \wedge t^{-2})
\,dt\nonumber\\[-8pt]\\[-8pt]
&\leq& C \varepsilon^{(\zeta- 1)^+}.\nonumber
\end{eqnarray}
It thus follows that
%
%
\begin{equation}\label{eqDeps-1}
\E\bigl\| \hD_{\varepsilon y} \bigl(\psitz(t)
- \psitz(s) \bigr) \bigr\|_{L^2}^2
\leq C |\varepsilon y|^{-1} \varepsilon^{(\zeta- 1)^+}.
\end{equation}
On the other hand, (\ref{eqpsi-1}) and (\ref{eqpsi-2}) imply that
\[
\E|\thtk(t) - \thtk(s)|^2 \leq C (1 \wedge k^{-2})^{3/4} |t -
s|^{1/4} ,
\]
and therefore
%
%
\begin{eqnarray}
\label{eqDeps-2}
\E\bigl\| \hD_{\varepsilon y} \bigl(\psitz(t)
- \psitz(s) \bigr) \bigr\|_{L^2}^2
& = & \sum_{|k| > \varepsilon^{-\zeta}}
\biggl|\frac{e^{ik\varepsilon y} -1 }{\varepsilon y}\biggr|^2
\E| \thtk(t) - \thtk(s) |^2
\nonumber\\
& \leq & C |\varepsilon y|^{-2} |t - s|^{1/4}
\sum_{|k| > \varepsilon^{-\zeta}} (1 \wedge k^{-2})^{3/4}
\\
& \leq & C |\varepsilon y|^{-2} |t - s|^{1/4} .\nonumber
\end{eqnarray}
Combining (\ref{eqDeps-1}) and (\ref{eqDeps-2}), we find that
\[
\E\bigl\| \hD_{\varepsilon y} \bigl(\psitz(t)
- \psitz(s) \bigr) \bigr\|_{L^2}^2
\leq
|\varepsilon y|^{-1- \kappa} |t - s|^{\kappa/4}
\varepsilon^{(1-\kappa)(\zeta- 1)^+}.
\]
Similarly, we obtain
\begin{eqnarray*}
\E\| \hD_{\varepsilon y} \psitz(0) \|_{L^2}^2
& = & \sum_{|k| > \varepsilon^{-\zeta}}
\biggl|\frac{e^{ik\varepsilon y} -1 }{\varepsilon y}\biggr|^2
\E| \thtk(0) |^2
\\
& \leq & C \sum_{|k| > \varepsilon^{-\zeta}}
\frac{1 - \cos(k\varepsilon y) }{|\varepsilon k y|^2}
\leq C |\varepsilon y|^{-1} \varepsilon^{(\zeta- 1)^+}.
\end{eqnarray*}
In view of Lemma~\ref{lemKolmogorov}, the latter two estimates imply that
\[
\E\sup_{t\in[0,T]}
\| \hD_{\varepsilon y} \psitz(t) \|_{L^2}^2
\leq C |\varepsilon y|^{-1 - \kappa} \varepsilon^{(1-\kappa)(\zeta- 1)^+}.
\]
Using this bound, the desired result for $\Theta(\psitz(t))$ can be
obtained easily, since
\begin{eqnarray*}
\E\sup_{t\in[0,T]} \Theta(\psitz(t))
& = &
\E\sup_{t\in[0,T]} \int_\R|y|^2
\| \hD_{\varepsilon y} \psitz(t) \|_{L^2}^2 |\mu|(dy)
\\& \leq &
\int_\R|y|^2 \E\sup_{t\in[0,T]}
\| \hD_{\varepsilon y} \psitz(t) \|_{L^2}^2 |\mu|(dy)
\\& \leq &
C \varepsilon^{-1 - \kappa+ (1-\kappa)(\zeta- 1)^+}
\int_\R|y|^{1 - \kappa} |\mu|(dy)
\\& \leq &
C \varepsilon^{-1 - \kappa+ (1-\kappa)(\zeta- 1)^+}.
\end{eqnarray*}
The result now follows by rescaling $\kappa$.
\end{pf}

The next and final result of this section involves the term which gives
rise to the correction term in the limiting equation.
Before stating the result, we introduce the notation
\begin{eqnarray*}
\Xi_\varepsilon^y(u) &\stackrel{\mathrm{def}}{=}& \frac{\varepsilon y^2}{2}
\hD_{\varepsilon y} u \ot\hD_{\varepsilon y} u ,
\\
\Lambda^y &\stackrel{\mathrm{def}}{=}& \frac{1}{2\pi\nu} \int_{\R_+}
\frac
{h^2(t)}{t^2 f(t)}
\bigl(1- \cos(yt)\bigr) \,dt
\end{eqnarray*}
and
\[
\Lambda_\varepsilon^y \stackrel{\mathrm{def}}{=}
\sum_{\varepsilon^{-\gamma} < k < \varepsilon^{-\chi}}
\Lambda_{\varepsilon,k}^y
\stackrel{\mathrm{def}}{=}
\sum_{\varepsilon^{-\gamma} < k < \varepsilon^{-\chi}} \frac{(1 -
\cos(\varepsilon k y))
h^2(\varepsilon k)}{2\pi\varepsilon(1 + \nu k^2f(\varepsilon k))}.
\]
Note that one has the identities
\[
\Xi_\varepsilon(u) = \int_\R\Xi_\varepsilon^y(u) \mu(dy) ,\qquad
\Lambda= \int_\R\Lambda^y \mu(dy) ,\qquad
\E\Xi_\varepsilon^y(\psitgx) = \Lambda_\varepsilon^y I ,
\]
where the constant $\Lambda$ is given by (\ref{edefLambda}).
\begin{proposition}\label{propsecond-chaos}
Let $\alpha> {1\over2}$, $\gamma\le{1\over2}$ and $\chi\ge
{3\over2}$. For $\varepsilon\in(0,1]$, we then have
\[
\E\sup_{t \in[0,T]} \| \Lambda- \Xi_\varepsilon(\psitgx(t)) \|
_{-\alpha}
\leq C \varepsilon^{1/2} .
\]
\end{proposition}
\begin{pf}
The proof is an application of Lemma~\ref{lemKolmogorov} with $\xi=
\Lambda- \Xi_\varepsilon(\psitgx)$. For brevity, we shall write $A
:= \Xi_\varepsilon(\psitgx)$ and $A^y := \Xi_\varepsilon^y(\psitgx
)$. We divide the proof into several steps.\vadjust{\goodbreak}

\textit{Step} 1.
First, we claim that $\xi(t) = \Lambda- A(t)$ satisfies the condition
(\ref{eqequiv-mom}) concerning the equivalence of all $q$-moments.

To see this, note that $\psitgx$ admits the representation $\psit(t)
= \sum_k \alpha_k(t) e_k$ where each $\alpha_k(t)$ is a~Gaussian
random vector in $\R^n$. As a~consequence, for every $y \in\R$, each
component of $\Lambda_\varepsilon^y - A^y$ is a~polynomial of
Gaussian random variables of degree at most two. It thus belongs to the
direct sum of Wiener chaoses of order $\leq2$ and the same is true for
$\Lambda_\varepsilon- A$, since each Wiener chaos is a~closed
subspace of the space of square integrable random variables. The claim
thus follows from the well-known equivalence of moments for Hilbert
space-valued Wiener chaos (see, e.g., \cite{KW92}).

\textit{Step} 2.
In this step, we estimate how well $\Lambda_\varepsilon^y$
approximates $\Lambda^y$.
Since $|1-\cos x| \le C(1\wedge x^2)$,
we have the bound
$|\Lambda_{\varepsilon,k}^y| \le C(\varepsilon y^2 \wedge
(\varepsilon k^2)^{-1})$ for some constant $C$. As an immediate
consequence,
we have the bound
%
%
\begin{equation}\label{eboundC}
\biggl|\Lambda_\varepsilon^y - \sum_{k \ge1}\Lambda_{\varepsilon
,k}^y\biggr| \le C (\varepsilon^{1-\gamma} y^2 + \varepsilon
^{\chi-1}) .
\end{equation}
Define now the function
\[
\Phi_y(t) = \frac{(1 - \cos(yt))
h^2(t)}{2\pi\nu t^2 f(t)} ,
\]
so that, since $h^2/f$ is bounded by assumption, we obtain the bound
\[
|\Lambda_{\varepsilon,k}^y - \varepsilon\Phi_y(\varepsilon k)| \le
C {\varepsilon y^2 \over k^2} .
\]
Combining this bound with (\ref{eboundC}), we have
\[
\biggl|\Lambda_\varepsilon^y - \sum_{k \ge1}\varepsilon\Phi
_y(\varepsilon k)\biggr| \le C (\varepsilon^{1-\gamma} y^2 +
\varepsilon^{\chi-1}) .
\]
At this stage, we recall that for any function $\Phi$ of bounded
variation, one has the approximation
\[
\biggl|\sum_{k \ge1}\varepsilon\Phi(\varepsilon k) - \int_0^\infty
\Phi(t) \,dt\biggr| \le\varepsilon\|\Phi\|_\BV,
\]
where $\|\Phi\|_\BV$ denotes the variation of $\Phi$ over $\R_+$.
Furthermore, for any pair~$\Phi$, $\Psi$, we have the bound
%
%
\begin{equation}\label{eprodBV}
\|\Phi\Psi\|_\BV\le\|\Phi\|_{L^\infty}\|\Psi\|_\BV+ \|\Psi\|
_{L^\infty}\|\Phi\|_\BV.
\end{equation}
If we set $\Psi_y(t) = (1 - \cos(yt))/t^2$, we have
\begin{eqnarray*}
\|\Psi_y\|_\BV &=& \int_0^\infty|\Psi_y'(t)| \,dt = \int_0^\infty
{|yt \sin yt + 2 \cos yt - 2| \over t^3} \,dt \\
&\le& C|y|^3\int_0^\infty\biggl(1 \wedge{1\over y^2 t^2}\biggr) \,dt
\le Cy^2 .
\end{eqnarray*}
Since $\Psi(0) = y^2/2$, a~similar bound holds for its $L^\infty$
norm, and we conclude from (\ref{eprodBV}) that
\[
\|\Phi_y\|_\BV\le C y^2 .
\]
It follows immediately that we have the bound
%
%
\begin{equation}\label{ediffCCeps}
|\Lambda_\varepsilon^y - \Lambda^y| \le C (\varepsilon
^{1-\gamma} y^2 + \varepsilon y^2 + \varepsilon^{\chi-1}) .
\end{equation}

\textit{Step} 3.
We now use these bounds in order to obtain control over $\|\Lambda- A
\|_{-\alpha}^2$ for a~fixed time $t \geq0$ (which is often suppressed
in the notation).

In order to shorten the notation, note
that, we can write
\[
\psit(x,t) = \sum_{k \in\Z} \frac{h(\varepsilon|k|)}{\sqrt
{2}\sqrt{1 +\nu k^2 f(\varepsilon|k|)}} \eta_k(t) e_k(x) ,
\]
where the $\eta_k$ are a~sequence of i.i.d. $\C^n$-valued
Ornstein--Uhlenbeck processes with
\[
\E\bigl(\eta_k(t) \otimes\eta_\ell(s)\bigr) = \CE_k^{t-s}
\delta_{k,-\ell}I , \qquad\CE_k^t = \exp\bigl(-\bigl(1+\nu k^2f(\varepsilon
|k|)\bigr)|t|\bigr) ,
\]
and satisfying the reality condition $\eta_{-k} = \bar\eta_k$. Here,
$I$ denotes the identity matrix.
We will also use the notational shortcut
\[
\mathbf{A}_{k,\ell}^t \stackrel{\mathrm{def}}{=}\eta_k(t) \otimes\eta_\ell
(t) .
\]

Set now
\[
q_\varepsilon^{k} = {e^{ik\varepsilon y}-1\over\sqrt{2}}
{h(\varepsilon|k|) \over\sqrt{1 +\nu k^2 f(\varepsilon|k|)}} ,
\]
as a~shorthand.
With all of this notation in place, it follows
from the definition of $\Lambda_\varepsilon^y$ that
\[
A^y(t)- \Lambda_\varepsilon^y I
= \sum_{\varepsilon^{-\gamma}< |k|,|\ell| \leq\varepsilon^{-\chi}}
q_\varepsilon^{k}q_\varepsilon^{\ell}
(\mathbf{A}_{k,\ell}^t - \delta_{k,-\ell} I) e_{k+\ell}.
\]
As a~consequence, we have the identity
\begin{eqnarray*}
&&\E\|\Lambda_\varepsilon^y I - A^y(t)\|_{-\alpha}^2\\
&&\qquad = \sum_{k \in\Z} (1+|k|^2)^{-\alpha}
\sum_{\ell, m} q_\varepsilon^{\ell}q_\varepsilon^{k-\ell}
\bar q_\varepsilon^{m}\bar q_\varepsilon^{k-m}
\\
&&\qquad\quad\hspace*{85pt}{}\times\E\tr\bigl((\mathbf{A}^t_{\ell,k-\ell} - \delta
_{k,0}I)
(\bar\mathbf{A}^t_{m,k-m} - \delta_{k,0} I)\bigr) ,
\end{eqnarray*}
where the second sum ranges over all $\ell, m \in\Z$ for which $\ell
, k-\ell, m, k-m$ belong to $(\varepsilon^{-\gamma}, \varepsilon
^{-\chi}]$.
A straightforward case analysis allows to check that
%
%
\begin{equation}\label{ecaseAnal}\qquad
\E\tr\bigl((\mathbf{A}^t_{\ell,k-\ell} - \delta_{k,0}I)
(\bar\mathbf{A}^t_{m,k-m} - \delta_{k,0} I)\bigr) = n\delta
_{\ell,m} + n^2 \delta_{\ell,k-m} ,
\end{equation}
so that
\[
\E\|\Lambda_\varepsilon^y I - A^y(t)\|_{-\alpha}^2 \leq C\sum_{k
\in\Z} (1+|k|^2)^{-\alpha}
\sum_{\ell\in\Z} |q_\varepsilon^{\ell}|^2 |q_\varepsilon^{k-\ell
}|^2 .\vadjust{\goodbreak}
\]
Note now that there exists a~constant $C$ such that the bound
\[
|q_\varepsilon^k| \le C\sqrt\varepsilon\biggl(|y| \wedge{1\over
\varepsilon|k|}\biggr) \le C \varepsilon^{(1-\beta)/2}
|k|^{-{\beta/2}} |y|^{1-{\beta/2}}
\]
is valid for all $\varepsilon< 1$, $k\in\Z$, $y \in\R$, and $\beta
\in[0,1]$. It follows that there exists a~constant $C>0$ such that we
have the bound
\begin{eqnarray*}
\E\|\Lambda_\varepsilon^y I - A^y(t)\|_{-\alpha}^2 &\le& C \sum
_{\ell, m \ge1} {|q_\varepsilon^\ell|^2|q_\varepsilon^m|^2\over
|\ell+ m|^{2\alpha}} \le
C \sum_{\ell, m \ge1} {|q_\varepsilon^\ell|^2|q_\varepsilon
^m|^2\over|\ell|^\alpha|m|^{\alpha}} \\
&\le& C \varepsilon^{2-2\beta}\sum_{\ell, m \ge1} {|y|^{4-2\beta}
\over|\ell|^{\alpha+ \beta}|m|^{\alpha+ \beta}} \le C \varepsilon
|y|^3 ,
\end{eqnarray*}
where we made the choice $\beta= {1\over2}$ to obtain the last bound,
using the fact that $\alpha> {1\over2}$ by assumption.
Combining this bound with (\ref{ediffCCeps}), the constraints $\gamma
\le{1\over2}$ and $\chi\ge{3\over2}$,
and using the fact that $\mu$ has finite fourth moment, we have
\begin{eqnarray*}
\E\|\Lambda I - A(t) \|_{-\alpha}^2 &\le& C \int\E\|\Lambda^y -
A^y(t)\|_{-\alpha}^2 |\mu|(dy) \\
&\le& C \int\E\|\Lambda^y_\varepsilon- A^y(t)\|_{-\alpha}^2 |\mu
|(dy) + C\varepsilon\\
&\le& C\varepsilon.
\end{eqnarray*}

\textit{Step} 4.
Finally, we shall estimate $\E\|A(t) - A(s)\|_{-\alpha}^2$.
Similarly to (\ref{ecaseAnal}), this involves the identity
\[
\E\tr( \mathbf{A}^t_{\ell,k-\ell}\bar\mathbf{A}^s_{m,k-m})
= n \delta_{k,0} + (n \delta_{l,m} + n^2 \delta_{l,k-m})
\CE_{\ell}^{t-s}\CE_{k-m}^{t-s} .
\]
As a~consequence, we infer that
\begin{eqnarray*}
D_{k\ell m}(t,s) &\stackrel{\mathrm{def}}{=}& \E\tr\bigl( (\mathbf{A}^t_{\ell
,k-\ell} -
\mathbf{A}^s_{\ell,k-\ell})
(\bar\mathbf{A}^t_{m,k-m} - \bar\mathbf{A}^s_{m,k-m}) \bigr)
\\& = & 2(n \delta_{\ell,m} + n^2 \delta_{\ell,k-m})
(1 - \CE_{\ell}^{t-s}\CE_{k-m}^{t-s}) .
\end{eqnarray*}
It thus follows that for any $\delta\in[0,1]$,
\begin{eqnarray*}
&& D_{k\ell m}(t,s)
\\
&&\qquad \leq C(\delta_{\ell,m} + \delta_{\ell,k-m})\bigl(1 \wedge
\bigl(2 + \nu\ell^2 f(\varepsilon|\ell|) + \nu(k-m)^2 f(\varepsilon
|k-m|)\bigr) |t-s|\bigr)
\\
&&\qquad \leq C(\delta_{\ell,m} + \delta_{\ell,k-m}) |t-s|^{\delta}
\bigl(1 + \ell^{2\delta} f(\varepsilon|\ell|)^\delta
+ (k-m)^{2\delta} f(\varepsilon|k-m|)^\delta
\bigr) .
\end{eqnarray*}
Using this bound, we obtain
\begin{eqnarray*}
\hspace*{-4pt}&&\E\|A^y(t) - A^y(s)\|_{-\alpha}^2 \\
\hspace*{-4pt}&&\qquad= \sum_{k \in\Z}
(1+|k|^2)^{-\alpha}
\sum_{\ell, m} q_\varepsilon^{\ell}q_\varepsilon^{k-\ell}
\bar q_\varepsilon^{m}\bar q_\varepsilon^{k-m}
D_{k\ell m}(t,s)
\\
\hspace*{-4pt}&&\qquad \leq \sum_{k \in\Z} (1+|k|^2)^{-\alpha}
\sum_{\ell} |q_\varepsilon^{\ell}|^2 |q_\varepsilon^{k-\ell}|^2
\bigl(D_{k\ell\ell}(t,s) + D_{k,\ell,k-\ell}(t,s)\bigr)
\\
\hspace*{-4pt}&&\qquad \leq C |t-s|^{\delta} \sum_{k \in\Z} (1+|k|^2)^{-\alpha}
\sum_{\ell} |q_\varepsilon^{\ell}|^2 |q_\varepsilon^{k-\ell}|^2
\\
\hspace*{-4pt}&&\qquad\quad\hspace*{127pt}{}
\times\bigl(1 + \ell^{2\delta} f(\varepsilon|\ell|)^\delta+
|k-\ell|^{2\delta}f(\varepsilon|k-\ell|)^\delta\bigr).
\end{eqnarray*}
Note that this expression is almost the same as in Step 3. Using the
calculations done there and taking into account that $h$ and $h/f^2$
are bounded functions, we infer that
\[
\E\|A^y(t) - A^y(s)\|_{-\alpha}^2
= C |t-s|^{\delta} \sum_{\ell, m \ge1} {|q_\varepsilon^\ell
|^2|q_\varepsilon^m|^2\over|\ell|^{\alpha- 2 \delta} |m|^{\alpha-
2 \delta}}
\le C |t-s|^{\delta} \varepsilon|y|^3 ,
\]
and therefore, using Jensen's inequality (which can be applied since
$|\mu|$ has finite mass), and Fubini's theorem,
\begin{eqnarray*}
\E\|A(t) - A(s)\|_{-\alpha}^2
& = &
\E\biggl\| \int_{\R} \bigl(A^y(t) - A^y(s) \bigr) \mu(dy) \biggr\|
_{-\alpha}^2
\\
& \leq & C
\int_{\R} \E\| A^y(t) - A^y(s)\|_{-\alpha}^2 |\mu|(dy)
\\
& \leq & C \varepsilon|t-s|^{\delta} \int_{\R} |y^3| |\mu|(dy)
\leq C \varepsilon|t-s|^{\delta} ,
\end{eqnarray*}
which is the desired bound.

The result follows by combining these steps with Lemma~\ref{lemKolmogorov}.
\end{pf}

\section{Convergence of the approximations}
\label{secapprox}

This last section is devoted to the convergence result itself. Recall that
we are considering a~number of intermediate processes $v_\varepsilon
^{(j)}$ with $j=1,\ldots, 4$
defined in (\ref{ev}). This section is correspondingly broken into
five subsections, with
the $j$th subsection yielding a~bound on $\|v_\varepsilon^{(j)} -
v_\varepsilon^{(j-1)}\|_\alpha$.
To prove these bounds, we shall introduce in each step a~stopping time
that forces the
difference between the processes considered in that step to remain
bounded. We then show that this difference actually
vanishes as $\varepsilon\to0$ with an explicit rate. As a~consequence,
the process actually does not ``see'' the stopping time with high probability.

\subsection{\texorpdfstring{From $\vbg$ to $v_\varepsilon^{(1)}$}{From v^gamma to v_epsilon^{(1)}}}

Define
\[
\tau_1^K := \tau^K \wedge
\inf\bigl\{ t \leq T \dvtx\bigl\|v_\varepsilon^{(1)}(t) - \vbg(t)\bigr\|_{\alpha}
\geq K \bigr\}.
\]
We shall show that for $t \leq\tau_K$, the $H^\alpha$-norm of
$v_\varepsilon^{(1)}(t) - \vbg(t)$ is controlled by the $L^\infty
$-norm of $\psiog$, which is of order $\varepsilon^{{\gamma/2}
- \kappa}$ for any $\kappa> 0$, as shown in Section \ref
{secprobtools}. The proof uses the mild formulations of the equations
for $v_\varepsilon^{(1)}$ and~$\vbg(t)$ as\vadjust{\goodbreak} well as the regularizing
properties of the semigroup $S$.
Note that\vspace*{1pt} the next proposition would still be true if we had replaced
the $H^\alpha$-norm in the definition of $\tau_1^K$ by the $L^\infty
$-norm. However, in the proof of Proposition~\ref{propest-12} below
it will be important to have a~bound on $v_\varepsilon^{(1)}$ in
$H^\alpha$.
\begin{proposition}\label{propest-01}
For $\kappa> 0$, we have
\[
\lim_{\varepsilon\to0} {\mathbb P}\Bigl( \sup_{t \leq\tau_1^K} \bigl\|
v_\varepsilon^{(1)}(t) - \vbg(t) \bigr\|_{\alpha} > \varepsilon^{
{\gamma/2} - \kappa} \Bigr) = 0 .
\]
\end{proposition}

\begin{pf}
Let $0 \leq s \leq t \leq\tau^*$. It follows from (\ref{eqv-0}) and
(\ref{eqv-1}) that $\varrho_\varepsilon:= v_\varepsilon^{(1)} -
\vbg$ satisfies the equation
\[
\varrho_\varepsilon(t) = S(t-s) \varrho_\varepsilon(s) + \int
_s^t S(t-r)
(\sigma^1 + \partial_x\sigma_\varepsilon^2)(r) \,dr ,
\]
where
\begin{eqnarray*}
\sigma_\varepsilon^1 & := & \tF(\vbg
+ \varrho_\varepsilon) - \tF( \vbg+ \psiog) ,\\
\sigma_\varepsilon^2 & := & G(\vbg
+ \varrho_\varepsilon) - G( \vbg+ \psiog) .
\end{eqnarray*}
Lemma~\ref{lemsemigroup-reg} yields the estimate
\begin{eqnarray*}
\|\varrho_\varepsilon(t)\|_{\alpha}
& \leq &\| \varrho_\varepsilon(s) \|_{\alpha}
+ C (t-s)^{(1-\alpha)/2}
\sup_{r\in(s,t)} \| ( \sigma_\varepsilon^1 + \partial_x \sigma
_\varepsilon^2)(r)\|_{-1}
\\
& \leq &\| \varrho_\varepsilon(s) \|_{\alpha}
+ C (t-s)^{(1-\alpha)/2}
\sup_{r\in(s,t)} \| \sigma_\varepsilon^1 (r)\|_{L^\infty}
+ \| \sigma_\varepsilon^2 (r)\|_{L^\infty} .
\end{eqnarray*}
Since $\vbg$, $\varrho_\varepsilon$, and $\psiog$ are bounded in
$L^\infty$-norm for $r \leq\tau_1^K$, and $F, G$ are $\CC^3$, it
follows that
\[
\| \sigma_\varepsilon^1 (r)\|_{L^\infty} + \| \sigma_\varepsilon^2
(r)\|_{L^\infty}
\leq C_K \| \varrho_\varepsilon(r) \|_{L^\infty} + C_K\| \psiog(r)
\|_{L^\infty},
\]
from which we infer that
\[
\|\varrho_\varepsilon(t)\|_{\alpha}
\leq\| \varrho_\varepsilon(s) \|_{\alpha}
+ C_K' (t-s)^{(1-\alpha)/2}
\sup_{r\in(s,t)}
\bigl( \| \varrho_\varepsilon(r) \|_{L^\infty}
+ \| \psiog(r) \|_{L^\infty}\bigr) .
\]
Choose $\delta_K > 0$ so small that $C_K' \delta_K^{(1-\alpha
)/2} \leq\frac12$, and set for $k \geq0$,
\[
r_k := \sup\{ \|\varrho_\varepsilon(t)\|_{\alpha}\dvtx
t \in[k \delta_K \wedge\tau_1^K,(k+1) \delta_K \wedge\tau_1^K]
\} .
\]
Taking into account that $H^\alpha\subseteq L^\infty$, we obtain the
inequality
\[
r_{k+1} \leq r_k + \frac12 r_{k+1}
+ \frac12 \sup_{t \in[0,T]} \| \psiog(t) \|_{L^\infty} ,
\]
which reduces to
\[
r_{k+1} \leq2 r_k
+ \sup_{t \in[0,T]} \| \psiog(t) \|_{L^\infty} .
\]
Combined with the estimate
\[
r_0 \leq2 \sup_{r \leq\delta_K \wedge\tau_1^K} \| \psiog(r) \|
_{L^\infty},\vadjust{\goodbreak}
\]
which can be derived similarly, it then follows that
\[
\sup_{t \in[0,\tau_1^K]} \| \varrho_\varepsilon(t)\|_\alpha
\leq\sup_{0 \leq k \leq T/\delta_K} r_k
\leq C_K \sup_{t \in[0,T]} \| \psiog(t) \|_{L^\infty} ,
\]
which together with Proposition~\ref{proppsi-est} implies the desired result.
\end{pf}

\subsection{\texorpdfstring{From $v_\varepsilon^{(1)}$ to $v_\varepsilon^{(2)}$}{From v_epsilon^{(1)} to v_epsilon^{(2)}}}

For the purpose of this section, we define the stopping time
\[
\tau_2^K := \tau_1^K \wedge
\inf\bigl\{ t \leq T \dvtx\bigl\|v_\varepsilon^{(2)}(t) - v_\varepsilon
^{(1)}(t)\bigr\|_{\alpha} \geq K
\bigr\}
\]
as well as the exponent
\[
\ta\stackrel{\mathrm{def}}{=}(1-\gamma) \alpha+ \frac{\gamma}{2} = {2\over
3} .
\]

\begin{proposition}\label{propest-12}
For $\kappa> 0$, we have
\[
\lim_{\varepsilon\to0} {\mathbb P}\Bigl( \sup_{t \leq\tau_2^K} \bigl\|
v_\varepsilon^{(2)}(t) - v_\varepsilon^{(1)}(t) \bigr\|_{\alpha} >
\varepsilon^{\ta- \kappa} \Bigr) = 0 .
\]
\end{proposition}

\begin{pf}
Let $0 \leq s \leq t \leq\tau^*$ and note that $\varrho_\varepsilon
:= v_\varepsilon^{(2)} - v_\varepsilon^{(1)}$ satisfies
\[
\varrho_\varepsilon(t) = S(t-s) \varrho_\varepsilon(s) + \int
_s^t S(t-r) \sigma_\varepsilon(r) \,dr ,
\]
where
\[
\sigma_\varepsilon:= \tF\bigl(v_\varepsilon^{(2)}\bigr) - \tF\bigl(v_\varepsilon^{(1)}\bigr)
+ D_\varepsilon\bigl(G\bigl(v_\varepsilon^{(1)} + \varrho_\varepsilon\bigr)\bigr)
- \partial_x G\bigl(v_\varepsilon^{(1)}\bigr) .
\]
From the definition of $\tau_2^K$, we know that $v_\varepsilon^{(1)}$
and $\varrho_\varepsilon$ are bounded in $L^\infty$ by a~constant
depending on $K$. Moreover, we have the bound $\| v_\varepsilon^{(1)}
\|_{\alpha} \leq C_K \varepsilon^{-\gamma(\alpha- 1/2) -
\kappa}$. Using these facts together with Corollary \ref
{corD-psi-eps-bound} we obtain, for $r \leq\tau_2^K$,
\begin{eqnarray*}
\| \sigma_\varepsilon\|_{-1}
& \leq &\bigl\| \tF\bigl(v_\varepsilon^{(2)}\bigr) - \tF\bigl(v_\varepsilon^{(1)}\bigr) \bigr\|
_{L^\infty}
\\
&&{} + \bigl\| (D_\varepsilon- \partial_x) G\bigl(v_\varepsilon^{(1)}\bigr) \bigr\|_{-1}
+ \bigl\| D^\varepsilon\bigl( G\bigl(v_\varepsilon^{(1)}
+ \varrho_\varepsilon\bigr) - G\bigl(v_\varepsilon^{(1)}\bigr) \bigr)\bigr\|_{-1}
\\
& \leq & C_K \|\varrho_\varepsilon\|_{L^\infty} + C_K \varepsilon
^\alpha\bigl( 1 + \bigl\| v_\varepsilon^{(1)} \bigr\|_{\alpha}\bigr)
+ \bigl\| G\bigl(v_\varepsilon^{(1)}
+ \varrho_\varepsilon\bigr) - G\bigl(v_\varepsilon^{(1)}\bigr) \bigr\|_{L^\infty}
\\
& \leq & C_K (\varepsilon^{\ta-\kappa}
+ \| \varrho_\varepsilon\|_{L^\infty}) ,
\end{eqnarray*}
hence
\begin{eqnarray*}
\| \varrho_\varepsilon(t) \|_\alpha
& \leq &\| \varrho_\varepsilon(s) \|_\alpha
+ C(t-s)^{(1-\alpha)/2}
\sup_{r \in(s,t)} \|\sigma_\varepsilon(r)\|_{-1}
\\
& \leq &\| \varrho_\varepsilon(s) \|_\alpha
+ C_K (t-s)^{(1-\alpha)/2}
\sup_{r \in(s,t)}
\bigl( \varepsilon^{\ta- \kappa}
+ \| \varrho_\varepsilon(r) \|_{L^\infty} \bigr) .\vadjust{\goodbreak}
\end{eqnarray*}
Arguing as in the proof of Proposition~\ref{propest-01}, it follows that
\[
\sup_{t\in[0,\tau_2^K]}
\| \varrho_\varepsilon(t) \|_\alpha
\leq C_K \varepsilon^{\ta- \kappa} ,\vspace*{-2pt}
\]
which immediately yields the desired result.\vspace*{-2pt}
\end{pf}

\subsection{\texorpdfstring{From $v_\varepsilon^{(2)}$ to $v_\varepsilon^{(3)}$}{From v_epsilon^{(2)} to v_epsilon^{(3)}}}

Define
\[
\tau_3^K := \tau_2^K \wedge
\inf\bigl\{ t \leq T\dvtx\bigl\|v_\varepsilon^{(3)}(t) - v_\varepsilon
^{(2)}(t)\bigr\|_{\alpha}
\geq K \bigr\}.\vspace*{-2pt}
\]

In this case, the singularity $(t-s)^{-\alpha/2}$ which arises
in the proof below,
prevents us from arguing as in Proposition~\ref{propest-01}. We
nevertheless have the following proposition.\vspace*{-2pt}
\begin{proposition}\label{propest-23}
For $\kappa> 0$, we have
\[
\lim_{\varepsilon\to0} {\mathbb P}\Bigl( \sup_{t \leq\tau_3^K} \bigl\|
v_\varepsilon^{(3)}(t) - v_\varepsilon^{(2)}(t) \bigr\|_{\alpha} >
\varepsilon^{\zeta- \kappa} \Bigr) = 0 ,\vspace*{-2pt}
\]
where the exponent $\zeta$ is given by
\[
\zeta\stackrel{\mathit{def}}{=}\ta\wedge\bigl(\tfrac32 - \alpha\bigr) \wedge\bigl(2 -
\gamma
\bigl(\alpha+ \tfrac32\bigr)\bigr) = \tfrac23 .\vspace*{-2pt}
\]
\end{proposition}
\begin{remark}
The exponent $\zeta$ arises by collecting the bounds (\ref
{eboundR1}), (\ref{eqR-2}), and (\ref{elastbound}).\vspace*{-2pt}
\end{remark}
\begin{pf*}{Proof of Proposition~\ref{propest-23}}
Let $0 \leq s \leq t \leq\tau^*$. It follows from (\ref{eqv-2}) and
(\ref{eqv-3}) that $\varrho_\varepsilon:= v_\varepsilon^{(3)} -
v_\varepsilon^{(2)}$ satisfies
\begin{eqnarray}
\label{eqSst}
\varrho_\varepsilon(t) &=& S_\varepsilon(t-s)\varrho_\varepsilon(s) +
\bigl( S_\varepsilon(t-s) - S(t-s) \bigr) v_\varepsilon^{(2)}(s)
\nonumber\\[-9pt]\\[-9pt]
&&{}+ \CR_1(s,t)
+ \CR_2(s,t) ,\nonumber\vspace*{-2pt}
\end{eqnarray}
where
\[
\CR_1(s,t) \stackrel{\mathrm{def}}{=}\bigl( \psitug(t) - \psioug(t) \bigr) - \bigl(
S_\varepsilon
(t-s)\psitug(s)
- S(t-s)\psioug(s) \bigr)\vspace*{-2pt}
\]
and
\begin{eqnarray*}
\CR_2(s,t)& := & \int_s^t \bigl( S_\varepsilon(t-r) - S(t-r) \bigr)
\bigl( \tF\bigl(v_\varepsilon^{(3)}(r)\bigr)
+ D_\varepsilon G\bigl( v_\varepsilon^{(3)}(r) \bigr)\bigr) \,dr
\\[-2pt]
&&{} + \int_s^t S(t-r) \bigl( \tF\bigl(v_\varepsilon^{(3)}(r)\bigr) -
\tF\bigl(v_\varepsilon^{(2)}(r)\bigr)
\\[-2pt]
&&\hspace*{67.7pt}{}
+ D_\varepsilon G\bigl( v_\varepsilon^{(3)}(r)\bigr)
- D_\varepsilon G\bigl( v_\varepsilon^{(2)}(r) \bigr) \bigr) \,dr .\vspace*{-2pt}
\end{eqnarray*}

We shall first prove a~bound on $\CR_1(s,t)$.
Using both inequalities from Lem\-ma~\ref{lemS-comparison}, we obtain
\begin{eqnarray*}
&& \bigl\| \bigl( S_\varepsilon(t-s)\psitug(s) - S(t-s)\psioug(s) \bigr) \bigr\|
_\alpha
\\[-2pt]
&&\qquad \leq
\bigl\| S_\varepsilon(t-s)
\bigl( \psitug(s) - \psioug(s) \bigr) \bigr\|_\alpha
+ \bigl\| \bigl( S_\varepsilon(t-s) - S(t-s)\bigr)\psioug(s) \bigr\|_\alpha
\\[-2pt]
&&\qquad \leq
C \bigl\|
\bigl( \psitug(s) - \psioug(s) \bigr) \bigr\|_\alpha
+ C \varepsilon^2 \|\psioug(s) \|_{\alpha+ 2} ,
\end{eqnarray*}
and therefore
\begin{eqnarray*}
\| \CR_1(s,t) \|_\alpha
& \leq &\bigl\|
\bigl( \psitug(t) - \psioug(t) \bigr) \bigr\|_\alpha
+ C \bigl\|
\bigl( \psitug(s) - \psioug(s) \bigr) \bigr\|_\alpha
\\[-2pt]
&&{}
+ C \varepsilon^2 \| \psioug(s) \|_{\alpha+2} .
\end{eqnarray*}
It thus follows from Proposition~\ref{proppsi-comparison} that
%
%
\begin{equation}\label{eboundR1}
\E\sup_{s, t\in[0,T]} \| \CR_1(s,t) \|_\alpha
\leq C \varepsilon^{2 - \gamma(\alpha+ 3/2) - \kappa} .
\end{equation}

We shall now prove a~bound on $\CR_2(s,t)$.
For this purpose, we note that the definitions of the various stopping
times imply that $v_\varepsilon^{(2)}(t)$ is bounded in $H^\alpha
$-norm by $C_K \varepsilon^{-\gamma(\alpha- 1/2)- \kappa}$.
Using this fact, together with Lemmas~\ref{lemS-comparison}, \ref
{lemD-eps-grad} and~\ref{lemsobolev-composition}, we obtain
\begin{eqnarray*}
&&
\biggl\| \int_s^t \bigl( S_\varepsilon(t-r) - S(t-r) \bigr)
\bigl( F\bigl( v_\varepsilon^{(3)}(r)\bigr)
+ D_\varepsilon G\bigl( v_\varepsilon^{(3)}(r) \bigr) \bigr) \,dr \biggr\|
_\alpha
\\
&&\qquad \leq\varepsilon^{\alpha}
\int_s^t (t-r)^{-(1+\alpha)/2}
\bigl\| F\bigl(v_\varepsilon^{(3)}(r)\bigr) + D_\varepsilon G\bigl( v_\varepsilon
^{(3)}(r) \bigr) \bigr\|_{\alpha-1} \,dr
\\
&&\qquad \leq C \varepsilon^{\alpha}
(t-s)^{(1-\alpha)/2}
\sup_{r \in[s,t]} \bigl(
\bigl\| F\bigl( v_\varepsilon^{(3)}(r) \bigr) \bigr\|_{\alpha} +
\bigl\| G\bigl( v_\varepsilon^{(3)}(r) \bigr) \bigr\|_{\alpha}\bigr)
\\
&&\qquad \leq C_K \varepsilon^{\alpha}
(t-s)^{(1-\alpha)/2}
\Bigl( 1+ \sup_{r \in[s,t]}
\bigl\| v_\varepsilon^{(3)}(r) \bigr\|_{\alpha} \Bigr)
\\
&&\qquad \leq C_K \varepsilon^{\alpha}
(t-s)^{(1-\alpha)/2}
\Bigl( 1+ \sup_{r \in[s,t]}
\bigl(\bigl\| v_\varepsilon^{(2)}(r) \bigr\|_{\alpha} + \| \varrho_\varepsilon(r)
\|_{\alpha}\bigr)\Bigr)
\\
&&\qquad \leq C_K \varepsilon^{\ta- \kappa}
(t-s)^{(1-\alpha)/2} .
\end{eqnarray*}
Furthermore, taking into account the $L^\infty$-bounds on
$v_\varepsilon^{(2)}$ and $\varrho_\varepsilon$ enforced by the
stopping times, Lemma~\ref{lemD-eps-grad} implies that
\begin{eqnarray*}
&&
\biggl\|\int_s^t  S(t-r)
\bigl(
F\bigl( v_\varepsilon^{(3)}(r) \bigr) - F\bigl( v_\varepsilon^{(2)}(r) \bigr)
+ D_\varepsilon G\bigl( v_\varepsilon^{(3)}(r) \bigr) - D_\varepsilon G\bigl(
v_\varepsilon^{(2)}(r) \bigr)
\bigr) \,dr \biggr\|_\alpha
\\
&&\qquad \leq C (t-s)^{(1-\alpha)/2}
\sup_{r\in(s,t)} \bigl\|
F\bigl( v_\varepsilon^{(3)}(r) \bigr) - F\bigl( v_\varepsilon^{(2)}(r) \bigr)
\\
&&\qquad\quad\hspace*{99pt}{} + D_\varepsilon G\bigl(
v_\varepsilon^{(3)}(r) \bigr)
- D_\varepsilon G\bigl( v_\varepsilon^{(2)}(r) \bigr) \bigr\|_{-1}
\\
&&\qquad\leq C (t-s)^{(1-\alpha)/2}
\sup_{r\in(s,t)} \bigl(\bigl\| F\bigl( v_\varepsilon^{(3)}(r) \bigr)
- F\bigl( v_\varepsilon^{(2)}(r) \bigr) \bigr\|_{L^\infty}
\\
&&\hspace*{98pt}\qquad\quad{}
+ \bigl\| G\bigl( v_\varepsilon^{(3)}(r) \bigr)
- G\bigl( v_\varepsilon^{(2)}(r) \bigr) \bigr\|_{L^\infty} \bigr)
\\
&&\qquad\leq C_K (t-s)^{(1-\alpha)/2}
\sup_{r\in(s,t)} \| \varrho_\varepsilon(r) \|_{L^\infty} .
\end{eqnarray*}
It thus follows that
%
%
\begin{equation} \label{eqR-2}
\|\CR_2(s,t) \|_\alpha\leq C_K' (t-s)^{(1-\alpha)/2}
\Bigl(\varepsilon^{\ta- \kappa} + \sup_{r\in(s,t)} \| \varrho
_\varepsilon(r) \|_{L^\infty} \Bigr) ,
\end{equation}
where we gave the constant a~name, since it will be reused below.

Choose $\delta_K \in(0,1)$ sufficiently small so that $C_K'
(\delta_K^{(1 - \alpha)/2} + \delta_K^{ \alpha/2}
) \leq\frac14$.
For $k \geq0$ put $\ell_k := k \delta_K \wedge\tau_3^K$,
and for $k \geq1$ set
\[
r_k := \sup\{ \|\varrho_\varepsilon(t)\|_{\alpha}\dvtx
t \in[ \ell_{k-1}, \ell_{k+1}] \} .
\]

Our next aim is to find a~bound for $r_1$. Observe that, when $s = 0$,
(\ref{eqSst}) simplifies to
%
%
\begin{equation}
\label{eqSst-simple}
\varrho_\varepsilon(t) = \bigl(S_\varepsilon(t) - S(t)\bigr) v_0
+\bigl( \psitug(t) - \psioug(t) \bigr)
+ \CR_2(0,t)
\end{equation}
with $\CR_2$ defined as previously.
Using Lemma~\ref{lemS-comparison} and the definition of $\tau^K$, we obtain
%
%
\begin{equation}\label{elastbound}
\bigl\| \bigl(S_\varepsilon(t) - S(t)\bigr) v_0 \bigr\|_\alpha
\leq C \varepsilon^{3/2 - \alpha- \kappa}
\| v_0 \|_{3/2 - \kappa}
\leq C_K \varepsilon^{3/2 - \alpha- \kappa} .
\end{equation}
Since $t \leq2 \delta_K$ and $C_K' \delta_K^{(1 - \alpha)/2}
\leq\frac14$, it follows from (\ref{eqR-2}) and (\ref
{eqSst-simple}) that
\[
r_1 \leq
C_K \varepsilon^{3/2 - \alpha- \kappa}
+ \sup_{t \in[0,T]}
\| \psitug(t) - \psioug(t) \|_\alpha
+ \frac12 ( \varepsilon^{\ta- \kappa} + r_1 ) ,
\]
hence, by definition of $\tau^K$,
%
%
\begin{equation} \label{eqr0-bound}
r_1 \leq C_K \varepsilon^{(3/2 - \alpha)
\wedge(2 - \gamma(\alpha+ 3/2))
\wedge\ta} \varepsilon^{-\kappa}
= C_K \varepsilon^{\zeta-\kappa} ,
\end{equation}
where $\zeta$ is defined as in the statement of the result.

Next, we shall prove a~recursive bound for $r_k$.
Note that the nonnegativity of the function $f$ in the definition of
$S_\varepsilon$ implies that
\[
\| S_\varepsilon(t-s)\varrho_\varepsilon(s) \|_\alpha
\leq\| \varrho_\varepsilon(s) \|_\alpha.
\]
Furthermore, by Lemma~\ref{lemS-comparison} and the fact that $\|
v_\varepsilon^{(2)}\|_\alpha\leq C_K \varepsilon^{-\gamma(\alpha
-1/2)- \kappa}$ before time $\tau_3^K$, we have
\begin{eqnarray*}
\bigl\| \bigl( S_\varepsilon(t-s) - S(t-s) \bigr) v_\varepsilon^{(2)}(s) \bigr\|_\alpha
& \leq & C \varepsilon^{\alpha} (t-s)^{-\alpha/2}
\bigl\| v_\varepsilon^{(2)}(s) \bigr\|_{\alpha}
\\
& \leq & C_K (t-s)^{-\alpha/2}\varepsilon^{\ta- \kappa} .
\end{eqnarray*}
Combining these bounds with (\ref{eqSst}) and (\ref{eqR-2}), we
find that
\begin{eqnarray*}
\|\varrho_\varepsilon(t) \|_\alpha
&\leq&\|\varrho_\varepsilon(s) \|_\alpha
+ C_K (t-s)^{- \alpha/2} \varepsilon^{\ta- \kappa}
+ \| \CR_1(s,t) \|_\alpha
\\
&&{}
+ C_K' (t-s)^{(1-\alpha)/2}
\Bigl( \varepsilon^{\ta- \kappa} + \sup_{r\in(s,t)} \| \varrho
_\varepsilon(r) \|_{\alpha} \Bigr) .
\end{eqnarray*}
Taking\vspace*{1pt} $k \geq1$, $s = \ell_{k-1}$, and $t \in[\ell_k, \ell
_{k+2}]$, it then follows, since $|t- s| \in[\delta_K, 3\delta_K]$
and $C_K'\delta^{(1-\alpha)/2} \leq\frac12$, that
\[
\| \varrho_\varepsilon(t)\|_\alpha
\leq\|\varrho_\varepsilon(\ell_{k-1}) \|_\alpha
+ C_K \varepsilon^{\ta- \kappa}
+ \| \CR_1(\ell_{k-1}, t) \|_\alpha
+ \tfrac12 \varepsilon^{\ta- \kappa} + \tfrac12 r_{k+1} .\vadjust{\goodbreak}
\]
Taking the supremum over $t \in[ \ell_k, \ell_{k+2}]$, we obtain
\[
r_{k+1}
\leq r_k
+ C_K \varepsilon^{\ta- \kappa}
+ \sup_{s,t \in[0,T]} \| \CR_1(s,t) \|_\alpha+ \frac12 r_{k+1} ,
\]
hence
%
%
\begin{equation} \label{eqrec}
r_{k+1}
\leq2 r_k
+ C_K \varepsilon^{\ta- \kappa}
+ 2 \sup_{s,t \in[0,T]} \| \CR_1(s,t) \|_\alpha.
\end{equation}

It readily follows from (\ref{eqr0-bound}) and (\ref{eqrec}) that
\[
\sup_{ t\in[0, \tau_3^K] } \| \varrho_\varepsilon(t) \|_\alpha
= \sup_{1\leq k \leq\lceil T/\delta_K \rceil} r_k
\leq C_K \Bigl( \varepsilon^{\zeta- \kappa}
+ \varepsilon^{\ta- \kappa}
+ \sup_{s,t \in[0,T]} \| \CR_1(s,t) \|_\alpha\Bigr) ,
\]
hence the result follows in view of the bound on $\CR_1(s,t)$.\vspace*{-3pt}
\end{pf*}

\subsection{\texorpdfstring{From $v_\varepsilon^{(3)}$ to $v_\varepsilon^{(4)}$}{From v_epsilon^{(3)} to v_epsilon^{(4)}}}

Define
\[
\tau_4^K := \tau_3^K \wedge
\inf\bigl\{ t \leq T\dvtx\bigl\|v_\varepsilon^{(4)}(t) - v_\varepsilon
^{(3)}(t)\bigr\|_{\alpha}
\geq K \bigr\} .\vspace*{-3pt}
\]

\begin{proposition}\label{propest-34}
For $\kappa> 0$, we have
\[
\lim_{\varepsilon\to0} {\mathbb P}\Bigl( \sup_{t \leq\tau_4^K} \bigl\|
v_\varepsilon^{(4)}(t) - v_\varepsilon^{(3)}(t) \bigr\|_{\alpha} >
\varepsilon^{\xi- \kappa} \Bigr) = 0 ,
\]
where
\[
\xi\stackrel{\mathit{def}}{=}\frac{\gamma}{2}
\wedge
\biggl(\ta- \frac12 \biggr)
\wedge
\biggl(\frac12 - \chi\biggl(\alpha- \frac12\biggr) \biggr) = {1\over8} .\vspace*{-3pt}
\]
\end{proposition}
\begin{remark}
Similarly to above, the exponent $\xi$ arises from the bounds (\ref
{eqsigma-1})--(\ref{eqsigma-4}).\vspace*{-3pt}
\end{remark}
\begin{pf*}{Proof of Proposition~\ref{propest-34}}
Let $0 \leq s \leq t \leq\tau^*$. It follows from (\ref{eqv-4}) and
(\ref{eqv-3}) that $\varrho_\varepsilon:= v_\varepsilon^{(4)} -
v_\varepsilon^{(3)}$ satisfies
\[
\varrho_\varepsilon(t) = S_\varepsilon(t-s)\varrho_\varepsilon(s)
+ \int_s^t S_\varepsilon(t-r) \sigma_\varepsilon(r)
\,dr ,
\]
where
\begin{eqnarray*}
\sigma_\varepsilon &:=& F\bigl(v_\varepsilon^{(4)} + \psitgx\bigr) -
F\bigl(v_\varepsilon^{(3)}\bigr)
\\
&&{} + \nabla G\bigl(v_\varepsilon^{(4)} + \psitgx\bigr)
D_\varepsilon\bigl(v_\varepsilon^{(4)} + \psitgx\bigr)
- D_\varepsilon G\bigl(v_\varepsilon^{(3)}\bigr)
+ \Lambda\Delta G\bigl(v_\varepsilon^{(3)}\bigr) .
\end{eqnarray*}
The definition of $D_\varepsilon$, together with (\ref{eqmoments}),
implies that for any function $u$ the following identity holds:
%
%
\begin{eqnarray}
\label{eqchain-approximate}
&&D_\varepsilon G(u)(x)\nonumber\\
&&\qquad = \nabla G(u(x)) D_\varepsilon u(x)
\nonumber\\
&&\qquad\quad{} + \int_\R\frac{\varepsilon y^2}{2}
D^2 G(u(x))[\hD_{\varepsilon y} u(x), \hD_{\varepsilon y} u(x)]
\mu(dy)
\\[-2pt]
&&\qquad\quad{} +
\int_\R\varepsilon^2 y^3 \int_0^{1} \int_0^t \int_0^s
D^3 G\bigl((1-r)u(x) + r u(x + \varepsilon y) \bigr)
\nonumber\\[-2pt]
&&\qquad\quad\hspace*{92pt}{}\times
[\hD_{\varepsilon y} u(x),\hD_{\varepsilon y} u(x),\hD_{\varepsilon y}u(x)]
\,dr \,ds \,dt\, \mu(dy) ,\nonumber
\end{eqnarray}
where the operator $\hD_\varepsilon$ is defined by taking $\mu:=
\delta_1 - \delta_0$ in the definition of~$D_\varepsilon$, that is,
$\hD_\varepsilon u(x) = \varepsilon^{-1}(u(x+\varepsilon) - u(x))$.
As a~consequence, we may write
\begin{eqnarray*}
\sigma_\varepsilon
& = & F\bigl(v_\varepsilon^{(4)} + \psitgx\bigr) - F\bigl(v_\varepsilon^{(3)}\bigr)
+ D_\varepsilon\bigl( G\bigl( v_\varepsilon^{(4)} + \psitgx\bigr) - G\bigl(v_\varepsilon^{(3)}\bigr)\bigr)
\\[-2pt]
&&{}
+ \bigl(\Lambda\Delta G\bigl(v_\varepsilon^{(3)}\bigr) - A\bigl(u_\varepsilon
^{(4)},u_\varepsilon^{(4)}\bigr) \bigr)
- B
\\[-2pt]
& = & F\bigl(v_\varepsilon^{(4)} + \psitgx\bigr) - F\bigl(v_\varepsilon^{(3)}\bigr)
+ D_\varepsilon\bigl( G\bigl( v_\varepsilon^{(4)} + \psitgx\bigr) - G\bigl(v_\varepsilon^{(3)}\bigr)\bigr)
\\[-2pt]
&&{}  - A\bigl(v_\varepsilon^{(4)},v_\varepsilon^{(4)}\bigr) - 2
A\bigl(v_\varepsilon^{(4)},\psitgx\bigr)
+ \bigl(\Lambda\Delta G\bigl(v_\varepsilon^{(3)}\bigr) - A(\psitgx,\psitgx) \bigr)
- B ,
\end{eqnarray*}
where we have used
\begin{eqnarray*}
A(v, w)(x) &\stackrel{\mathrm{def}}{=}&
\int_\R\frac{\varepsilon y^2}{2}
D^2 G \bigl(u_\varepsilon^{(4)}(x)\bigr)
[ \hD_{\varepsilon y} v(x), \hD_{\varepsilon y} w(x) ] \mu(dy)
,\\[-2pt]
B(x) &\stackrel{\mathrm{def}}{=}& \int_\R\varepsilon^2 y^3 \int_0^{1} \int
_0^t \int_0^s
D^3 G\bigl((1-r)u_\varepsilon^{(4)}(x) + r u_\varepsilon^{(4)}(x +
\varepsilon y) \bigr)
\\[-2pt]
&&\hspace*{82.4pt}{}\times
\bigl[\hD_{\varepsilon y} u_\varepsilon^{(4)}(x),\hD_{\varepsilon y}
u_\varepsilon^{(4)}(x),\\[-2pt]
&&\hspace*{157pt}\hD_{\varepsilon y}u_\varepsilon^{(4)}(x)\bigr]
\,dr \,ds \,dt \,\mu(dy)
\end{eqnarray*}
and $u_\varepsilon^{(4)} := v_\varepsilon^{(4)} + \psitgx$.

Our next aim is to prove the estimates (\ref{eqsigma-1})--(\ref
{eqsigma-4}) below in order to bound $\| \sigma_\varepsilon\|_{-1}$.

\subsubsection*{First term} Since $v_\varepsilon^{(4)}$, $\psitgx$
and $\varrho_\varepsilon$ are bounded in $L^\infty$ by definition of~$\tau_4^K$, it follows that
%
%
\begin{eqnarray}\label{eqsigma-1}
\bigl\| F\bigl( v_\varepsilon^{(4)} + \psitgx\bigr)
- F\bigl(v_\varepsilon^{(3)}\bigr) \bigr\|_{-1}
& \leq &C \bigl\| F\bigl( v_\varepsilon^{(4)} + \psitgx\bigr)
- F\bigl(v_\varepsilon^{(4)} - \varrho_\varepsilon\bigr) \bigr\|_{L^\infty}
\nonumber\\[-2pt]
& \leq &C_K ( \| \psitgx\|_{L^\infty}
+ \| \varrho_\varepsilon\|_{L^\infty})
\\[-2pt]
& \leq &C_K ( \varepsilon^{{\gamma/2} - \kappa}
+ \| \varrho_\varepsilon\|_{\alpha}) .\nonumber
\end{eqnarray}

\subsubsection*{Second term} We use Lemma~\ref{lemD-eps-grad} and
the fact that $v_\varepsilon^{(4)}$, $\psitgx$ and $\varrho
_\varepsilon$ are bounded in $L^\infty$ by definition of $\tau_4^K$,
to estimate
%
%
\begin{eqnarray}\label{eqsigma-t}\qquad
\bigl\| D_\varepsilon\bigl( G\bigl( v_\varepsilon^{(4)} + \psitgx\bigr)
- G\bigl(v_\varepsilon^{(3)}\bigr)\bigr) \bigr\|_{-1}
& \leq & C \bigl\| G\bigl( v_\varepsilon^{(4)} + \psitgx\bigr)
- G\bigl(v_\varepsilon^{(4)} - \varrho_\varepsilon\bigr) \bigr\|_{L^\infty}
\nonumber\\[-2pt]
& \leq & C_K ( \| \psitgx\|_{L^\infty}
+ \| \varrho_\varepsilon\|_{L^\infty})
\\[-2pt]
& \leq & C_K ( \varepsilon^{{\gamma/2} - \kappa}
+ \| \varrho_\varepsilon\|_{\alpha}) .\nonumber
\end{eqnarray}

\subsubsection*{Third and fourth term}
First, we note that for arbitrary functions $v, w$, one has
\[
\|A(v,w)\|_{-1}
\leq C \| A(v,w)\|_{L^1}
\leq C \varepsilon\bigl\|D^2 G\bigl(u_\varepsilon^{(4)}\bigr)\bigr\|_{L^\infty} \sqrt
{\Theta_\varepsilon(v) \Theta_\varepsilon(w)} .
\]
Since $\| u_\varepsilon^{(4)}(t) \|_{L^\infty} \leq C_K$ for $t \leq
\tau_4^K$, we have
\[
\bigl\| D^2 G\bigl(u_\varepsilon^{(4)}\bigr) \bigr\|_{L^\infty} \leq
C_K .
\]
Furthermore, we observe that $\|v_\varepsilon^{(4)}\|_\alpha\leq C_K
\varepsilon^{-\gamma(\alpha- 1/2) - \kappa}$ before time $\tau_4^K$.
Using this bound together with Lemma~\ref{lemD-eps-grad} and (\ref
{eqmoments}), we estimate
\begin{eqnarray*}
\Theta_\varepsilon\bigl(v_\varepsilon^{(4)}\bigr)
& = &\int_\R y^2 \bigl\| \hD_{\varepsilon y} v_\varepsilon^{(4)} \bigr\|
_{L^2}^2 |\mu|(dy)
\\& \leq &
C_K
\int_\R y^2 |\varepsilon y|^{2(\alpha-1)}
\bigl\| v_\varepsilon^{(4)} \bigr\|_\alpha^2 |\mu|(dy)
\\& \leq &
C_K
\int_\R y^2 |\varepsilon y|^{2(\alpha-1)}
\varepsilon^{-2\gamma(\alpha- 1/2)-2\kappa} |\mu|(dy)
\\& \leq & C_K \varepsilon^{2\ta- 2 - \kappa} .
\end{eqnarray*}
Moreover, by definition of the stopping time $\tau^K$ we have
\[
\Theta_\varepsilon(\psitgx) \leq\Theta_\varepsilon(\psitg) \leq
C_K \varepsilon^{-1 - \kappa}.
\]
Putting everything together, we obtain
%
%
\begin{equation} \label{eqsigma-2}
\bigl\| A\bigl(v_\varepsilon^{(4)}, v_\varepsilon^{(4)}\bigr) \bigr\|_{-1}
\leq C_K \varepsilon^{2\ta- 1 - 2 \kappa}
\end{equation}
and
%
%
\begin{equation} \label{eqsigma-3}
\bigl\| A\bigl(v_\varepsilon^{(4)}, \psitgx\bigr) \bigr\|_{-1} \leq C_K \varepsilon^{\ta
- 1/2 - 2\kappa} .
\end{equation}

\subsubsection*{Fifth term} Finally, we estimate $\Lambda\Delta G(
v_\varepsilon^{(3)}) - A(\psitgx, \psitgx) $.
By definition of $\tau_\varepsilon^K$, we have
$
\| \psitgx\|_\alpha\leq C_K \varepsilon^{-\chi(\alpha- 1/2)
- \kappa}
$
before time $\tau_4^K$. Since
$
\|v_\varepsilon^{(4)} \|_\alpha\leq C_K \varepsilon^{-\gamma(\alpha
- 1/2) - \kappa}
$
as observed before, we thus have
\[
\bigl\| u_\varepsilon^{(4)} \bigr\|_\alpha\leq C_K \varepsilon^{-\chi(\alpha
- 1/2) - \kappa} .
\]
Furthermore, since $\alpha> {1\over2}$, there exists
a constant $C>0$ such that we have the bound
\begin{eqnarray*}
\bigl\| \Lambda\Delta G\bigl(u_\varepsilon^{(4)}\bigr) - A(\psitgx, \psitgx)\bigr\|
_{-\alpha}
&=& \bigl\| \tr\bigl(D^2 G\bigl(u_\varepsilon^{(4)}\bigr)
\bigl(\Lambda I - \Xi_\varepsilon(\psitgx) \bigr)
\bigr) \bigr\|_{-\alpha}
\\& \le & C \bigl\|D^2 G\bigl(u_\varepsilon^{(4)}\bigr)\bigr\|_\alpha
\| \Lambda I - \Xi_\varepsilon(\psitgx) \|_{-\alpha} .
\end{eqnarray*}
Since the stopping time $\tau^K$ enforces that $ \| \Lambda I - \Xi
_\varepsilon(\psitgx) \|_{-\alpha} \leq C_K \varepsilon^{
1/2}$, we infer that
\[
\bigl\| \Lambda\Delta G\bigl(u_\varepsilon^{(4)}\bigr) - A(\psitgx, \psitgx)\bigr\|
_{-\alpha}
\le C_K \varepsilon^{1/2-\chi(\alpha- 1/2) - \kappa} .
\]
Since $u_\varepsilon^{(4)} - v_\varepsilon^{(3)} = \varrho
_\varepsilon+ \psitgx$, we have by definition of $\tau^K$,
\begin{eqnarray*}
\bigl\| \Delta G(u_\varepsilon^{4}) - \Delta G\bigl( v_\varepsilon^{(3)}\bigr)\bigr\|
_{-\alpha}
&\leq&\bigl\| \Delta G(u_\varepsilon^{4}) - \Delta G\bigl( v_\varepsilon
^{(3)}\bigr)\bigr\|_{L^\infty}
\\& \leq &C_K \| \varrho_\varepsilon\|_{L^\infty} + C_K\| \psitgx\|
_{L^\infty}
\\& \leq &C_K \| \varrho_\varepsilon\|_{\alpha} + C_K \varepsilon
^{{\gamma/2} - \kappa} .
\end{eqnarray*}
Putting these bounds together, we obtain
\begin{eqnarray}
\label{eqsigma-5}
&&\bigl\| \Lambda\Delta G\bigl(v_\varepsilon^{(3)}\bigr) - A(\psitgx, \psitgx)\bigr\|
_{-\alpha}\nonumber\\[-8pt]\\[-8pt]
&&\qquad\le C_K \bigl(\varepsilon^{1/2-\chi(\alpha- 1/2) - \kappa}
+ \varepsilon^{{\gamma}/{2} - \kappa}
+ \| \varrho_\varepsilon\|_{\alpha} \bigr) .\nonumber
\end{eqnarray}

\subsubsection*{Sixth term} To estimate $B$, we use the fact that $\|
u_\varepsilon^{(4)}(t)\|_{L^\infty} \leq C_K$ for $t \leq\tau_4^K$,
so that one has the bound
\[
\| B \|_{L^1}
\leq C_K
\iint_\R\varepsilon^2 y^3 \bigl| \hD_{\varepsilon y} u_\varepsilon
^{(4)}(x)\bigr|^3 |\mu|(dy) \, dx.
\]
We will split this expression into two parts, using the fact that
$u_\varepsilon^{(4)} = \psitgx+ v_\varepsilon^{(4)}$. First, using
the fact that $\Theta(\psitgx) \leq C_K \varepsilon^{1-\kappa}$
before time $\tau_4^K$ by definition of the stopping time $\tau^K$,
we find that
\begin{eqnarray*}
&& \iint_\R\varepsilon^2 y^3 |\hD_{\varepsilon y} \psitgx(x)|^3
|\mu|(dy) \, dx
\\
&&\qquad
\leq2 \| \psitgx\|_{L^\infty}
\iint_\R\varepsilon y^2 | \hD_{\varepsilon y}\psitgx(x)|^2
|\mu|(dy) \, dx
\\
&&\qquad = 2 \varepsilon\| \psitgx\|_{L^\infty}
\Theta(\psitgx)
\leq C_K \varepsilon^{-2\kappa} \| \psitgx\|_{L^\infty}
\leq C_K \varepsilon^{{\gamma}/{2} - 3\kappa} .
\end{eqnarray*}
Second,\vspace*{1pt} using the fact that $H^{1/6} \subseteq L^3$, Lemma \ref
{lemD-eps-grad} and the fact that $\|v_\varepsilon^{(4)} \|_\alpha
\leq C_K \varepsilon^{-\gamma(\alpha- 1/2) - \kappa}$, we obtain
\begin{eqnarray*}
&&\iint_\R\varepsilon^2 y^3 \bigl| \hD_{\varepsilon y} v_\varepsilon
^{(4)}(x)\bigr|^3 |\mu|(dy) \, dx\\
&&\qquad \leq  C \varepsilon^2 \int_\R|y|^3 \bigl\| \hD_{\varepsilon y}
v_\varepsilon^{(4)} \bigr\|_{1/6}^3 |\mu|(dy)
\\
&&\qquad\leq C\varepsilon^{3\alpha- 3/2} \int_\R|y|^{3 \alpha-
1/2} |\mu|(dy)
\bigl\| v_\varepsilon^{(4)} \bigr\|_{\alpha}^3
\\
&&\qquad\leq C_K \varepsilon^{3\ta- 3/2 - 3 \kappa} .
\end{eqnarray*}
It thus follows that
%
%
\begin{equation} \label{eqsigma-4}
\| B \|_{L^1}
\leq C_K \varepsilon^{ {\gamma}/{2} - 3\kappa}
+ C_K \varepsilon^{3(\ta- 1/2) - 3 \kappa} .
\end{equation}

Combining inequalities (\ref{eqsigma-1})--(\ref{eqsigma-5}), we find
that
\begin{eqnarray*}
\sup_{r \in(s,t)} \| \sigma_\varepsilon(r)\|_{-1}
&\leq& C_K \Bigl(
\varepsilon^{{\gamma}/{2} - 3 \kappa}
+ \varepsilon^{\ta- 1/2 - 2\kappa}\\[-2pt]
&&\hspace*{18.6pt}{}
+ \varepsilon^{1/2-\chi(\alpha- 1/2) - \kappa}
+\sup_{r \in(s,t)} \| \varrho_\varepsilon(r) \|_{\alpha}
\Bigr) ,\vspace*{-2pt}
\end{eqnarray*}
and the result now follows as in Proposition~\ref{propest-12}.\vspace*{-3pt}
\end{pf*}

\subsection{\texorpdfstring{From $v_\varepsilon^{(4)}$ to $\vtg$}{From v_epsilon^{(4)} to v^gamma}}

Define
\[
\tau_5^K := \tau_4^K \wedge
\inf\bigl\{ t \leq T \dvtx\bigl\| \vtg(t) - v_\varepsilon^{(4)}(t)\bigr\|_{\alpha}
\geq K \bigr\}.\vspace*{-3pt}
\]

\begin{proposition}\label{propest-4eps}
For $\kappa> 0$, we have
\[
\lim_{\varepsilon\to0} {\mathbb P}\Bigl( \sup_{t \leq\tau_5^K} \bigl\|
\vtg(t) - v_\varepsilon^{(4)}(t) \bigr\|_{\alpha} > \varepsilon^{
1/2 \chi- 1/2 - \kappa} \Bigr) = 0 .\vspace*{-3pt}
\]
\end{proposition}
\begin{pf}
Let $0 \leq s \leq t \leq\tau^*$. It follows from (\ref{eqv-eps})
and (\ref{eqv-4}) that $\varrho_\varepsilon:= \vtg- v_\varepsilon
^{(4)}$ satisfies
\[
\varrho_\varepsilon(t) =
\int_s^t S_\varepsilon(t-r) \sigma_\varepsilon(r)
\,dr ,\vspace*{-2pt}
\]
where
\begin{eqnarray*}
\sigma_\varepsilon
& := &\nabla G( \vtg+ \psitg) D_\varepsilon(\vtg+ \psitg)\\[-2pt]
&&{} - \nabla G( \vtg+ \psitgx- \varrho_\varepsilon)
D_\varepsilon(\vtg+ \psitgx- \varrho_\varepsilon)
\\[-2pt]
&&{} + F(\vtg+ \psitg)
- F(\vtg+ \psitgx- \varrho_\varepsilon) .\vspace*{-2pt}
\end{eqnarray*}
In order to estimate $\sigma_\varepsilon$, we use (\ref
{eqchain-approximate}) to write
\begin{eqnarray*}
\sigma_\varepsilon
& = & D_\varepsilon G(\vtg+ \psitg)
- D_\varepsilon G(\vtg+ \psitgx- \varrho_\varepsilon)
\\[-2pt]
&&{} - \int_\R\frac{\varepsilon y^2}{2}
\bigl(D^2 G(u_\varepsilon) - D^2 G\bigl(u_\varepsilon^{(4)}\bigr) \bigr)
[\hD_{\varepsilon y} u_\varepsilon, \hD_{\varepsilon y}
u_\varepsilon] \mu(dy)
\\[-2pt]
&&{} - \int_\R\frac{\varepsilon y^2}{2}
D^2 G\bigl(u_\varepsilon^{(4)}\bigr)
\bigl[\hD_{\varepsilon y} \bigl(u_\varepsilon+ u_\varepsilon^{(4)}\bigr) , \hD
_{\varepsilon y} \bigl(u_\varepsilon- u_\varepsilon^{(4)}\bigr)\bigr] \mu(dy)
\\[-2pt]
&&{} - \varepsilon^2 \bigl( R_\varepsilon(u_\varepsilon
,u_\varepsilon) - R_\varepsilon\bigl(u_\varepsilon^{(4)},u_\varepsilon
^{(4)}\bigr) \bigr)
\\[-2pt]
&&{} + F(\vtg+ \psitg)
- F(\vtg+ \psitgx- \varrho_\varepsilon)
\\[-2pt]
& =&\!: \sigma_{\varepsilon,1} + \cdots+ \sigma_{\varepsilon,5} ,\vspace*{-2pt}
\end{eqnarray*}
where $u_\varepsilon:= \vtg+ \psitg$, $u_\varepsilon^{(4)} :=
v_\varepsilon^{(4)} + \psitgx$ and
\begin{eqnarray*}
R_\varepsilon(u^1, u^2)(x) &:=& \int_\R\varepsilon^2 y^3 \int_0^1
\int_0^t
\int_0^s D^3 G\bigl((1-r)u^1(x) +r u^1(x+\varepsilon y)\bigr)
\\[-2pt]
&&\hspace*{79.6pt}{}\times [\hD_{\varepsilon y} u^2, \hD
_{\varepsilon y} u^2,\hD_{\varepsilon y} u^2]
\,dr \,ds \,dt \,d \mu(y) .
\end{eqnarray*}
We shall now estimate $\sigma_{\varepsilon,1}, \ldots, \sigma
_{\varepsilon,5}$ individually.\vadjust{\goodbreak}

\subsubsection*{First term} First, we observe that $\vtg$, $\psitg
$, $\psitgx$ and $\varrho_\varepsilon$ are bounded in~$L^\infty$
before time $\tau_5^K$. Using Lemma~\ref{lemD-eps-grad}, the
embedding $H^\alpha\subseteq L^\infty$, and the definition of the
stopping time to bound $\| \psitlx\|_{L^\infty}$, we obtain
\begin{eqnarray}
\label{eqv4-1}
\| \sigma_{\varepsilon,1} \|_{-1}
& = &\bigl\| D_\varepsilon\bigl( G( \vtg+ \psitg)
- G(\vtg+ \psitgx- \varrho_\varepsilon)\bigr) \bigr\|_{-1}
\nonumber\\
& \leq & C \| G( \vtg+ \psitg)
- G(\vtg+ \psitgx- \varrho_\varepsilon) \|_{L^\infty}
\nonumber\\[-8pt]\\[-8pt]
& \leq & C_K ( \| \psitlx\|_{L^\infty}
+ \| \varrho_\varepsilon\|_{L^\infty})
\nonumber\\
& \leq & C_K ( \varepsilon^{\chi/2- \kappa}
+ \| \varrho_\varepsilon\|_{\alpha}) .\nonumber
\end{eqnarray}

\subsubsection*{Second term}
Using Lemma~\ref{lemD-eps-grad} and the fact that $\varepsilon
^{\gamma(\alpha- 1/2)+\kappa}\|\vtg\|_\alpha$ is bounded
before time $\tau_5^K$, we estimate
\begin{eqnarray*}
\Theta_\varepsilon(\vtg)
& = &\int_\R y^2 \| \hD_{\varepsilon y} \vtg\|_{L^2}^2 |\mu|(dy)
\\& \leq &
C_K
\int_\R y^2 |\varepsilon y|^{2(\alpha-1)}
\| \vtg\|_\alpha^2 |\mu|(dy)
\\& \leq &
C_K
\int_\R y^2 |\varepsilon y|^{2(\alpha-1)}
\varepsilon^{-2\gamma(\alpha- 1/2)-2\kappa} |\mu|(dy)
\\& \leq & C_K \varepsilon^{2\ta- 2 - \kappa} ,
\end{eqnarray*}
and by the definition of the stopping time $\tau^K$,
\[
\Theta_\varepsilon(\psitg) \leq C_K \varepsilon^{-1 - \kappa} .
\]
As a~consequence,
%
%
\begin{equation}
\label{eqTheta}\qquad
\Theta_\varepsilon(u_\varepsilon)
\leq
2\bigl( \Theta_\varepsilon(\vtg)
+ \Theta_\varepsilon(\psitg)\bigr)
\leq C_K ( \varepsilon^{2\ta- 2 - \kappa} + \varepsilon^{-1 -
\kappa} )
\leq\varepsilon^{-1 - \kappa} .
\end{equation}
Note that $\|u_\varepsilon\|_{L^\infty}$ and $\|u_\varepsilon^{(4)}\|
_{L^\infty}$ are bounded before time $\tau_5^K$.
Using that $L^1 \subseteq H^{-1}$, we obtain
\begin{eqnarray}
\label{eqv4-2}
\| \sigma_{\varepsilon,2} \|_{-1}
&\leq&\| \sigma_{\varepsilon,2} \|_{L^1}
\leq\varepsilon
\bigl\| D^2 G(u_\varepsilon) - D^2G\bigl(u_\varepsilon^{(4)}\bigr) \bigr\|_{L^\infty}
\Theta_\varepsilon(u_\varepsilon)
\nonumber\\
& \leq & C_K \varepsilon^{-2\kappa} \bigl\| u_\varepsilon- u_\varepsilon
^{(4)} \bigr\|_{L^\infty}
\nonumber\\[-8pt]\\[-8pt]
& \leq & C_K \varepsilon^{-2\kappa}
( \| \varrho_\varepsilon\|_{L^\infty}
+ \| \psitlx\|_{L^\infty}
)
\nonumber\\
& \leq & C_K \varepsilon^{-2\kappa}
( \| \varrho_\varepsilon\|_{\alpha}
+ \varepsilon^{\chi/2- \kappa}
) .
\nonumber
\end{eqnarray}

\subsubsection*{Third term}
By Lemma~\ref{lemD-eps-grad}, we have
\begin{eqnarray}
\label{eqTheta-rho}
\Theta_\varepsilon(\varrho_\varepsilon)
& = &\int_\R y^2 \| \hD_{\varepsilon y} \varrho_\varepsilon\|
_{L^2}^2 |\mu|(dy)
\nonumber\\[-8pt]\\[-8pt]
& \leq & C \int_\R y^2 |\varepsilon y|^{2(\alpha-1)}
\|\varrho_\varepsilon\|_\alpha^2 |\mu|(dy)
\leq C \varepsilon^{2\alpha-2}\|\varrho_\varepsilon\|_\alpha^2 .
\nonumber
\end{eqnarray}
Observe that $u_\varepsilon+ u_\varepsilon^{(4)} = 2 \vtg- \varrho
_\varepsilon+ \psitg+ \psitgx$ and
$u_\varepsilon- u_\varepsilon^{(4)} = \psitlx+ \varrho_\varepsilon$.
Taking into account that
\[
\varepsilon^{1 + \kappa} \Theta_\varepsilon(\psitg),\qquad
\varepsilon^{1+ \kappa} \Theta_\varepsilon(\psitg),\qquad
\varepsilon^{2 - \chi+ \kappa} \Theta_\varepsilon(\psitlx),\qquad
\| \varrho_\varepsilon\|_\alpha
\]
are all bounded before time $\tau_5^K$, we obtain
\begin{eqnarray}
\label{eqTheta-sum}\qquad
\Theta_\varepsilon\bigl(u_\varepsilon+ u_\varepsilon^{(4)}\bigr)
&\leq& C\bigl( \Theta_\varepsilon(\vtg) + \Theta_\varepsilon
(\varrho_\varepsilon)
+ \Theta_\varepsilon(\psitg) + \Theta_\varepsilon(\psitgx)\bigr)
\nonumber\\[-8pt]\\[-8pt]
& \leq& C_K (
\varepsilon^{2\ta- 2 - \kappa}
+ \varepsilon^{2\alpha- 1}
+ \varepsilon^{-1 - \kappa} )
\leq C_K \varepsilon^{-1 - \kappa}
\nonumber
\end{eqnarray}
and
\begin{eqnarray}
\label{eqTheta-diff}
\Theta_\varepsilon\bigl(u_\varepsilon- u_\varepsilon^{(4)}\bigr)
&\leq& C\bigl(\Theta_\varepsilon(\psitlx) + \Theta_\varepsilon
(\varrho_\varepsilon)\bigr)
\nonumber\\[-8pt]\\[-8pt]
&\leq& C_K \varepsilon^{\chi- 2 - \kappa}
+ C \varepsilon^{2\alpha-2} \|\varrho_\varepsilon\|_\alpha^2 .
\nonumber
\end{eqnarray}
Using that $\|u_\varepsilon^{(4)} \|_{L^\infty}\leq C_K$ before time
$\tau_5^K$, we obtain
\begin{eqnarray*}
\| \sigma_{\varepsilon,3} \|_{-1}
&\leq&
\| \sigma_{\varepsilon,3} \|_{L^1}
\leq\varepsilon
\bigl\| D^2G\bigl(u_\varepsilon^{(4)}\bigr) \bigr\|_{L^\infty}
\sqrt{ \Theta_\varepsilon\bigl(u_\varepsilon+ u_\varepsilon^{(4)}\bigr)
\Theta_\varepsilon\bigl(u_\varepsilon- u_\varepsilon^{(4)}\bigr)}
\\& \leq& C_K ( \varepsilon^{\chi/2- 1/2 - 2 \kappa}
+ \varepsilon^{\alpha- 1/2 - \kappa} \|\varrho_\varepsilon\|
_\alpha) .
\end{eqnarray*}

\subsubsection*{Fourth term}
We shall show that
%
%
\begin{equation}
\label{eqv4-4}
\| \sigma_{\varepsilon,4} \|_{-1}
\leq C_K (
\varepsilon^{-2\kappa}\| \varrho_\varepsilon\|_{\alpha}
+\varepsilon^{\chi/2- 1/2 - 2 \kappa}
).
\end{equation}

First, we use the $L^\infty$-bound on $u_\varepsilon^{(4)}$ enforced
by the stopping time, to obtain the pointwise bound
\begin{eqnarray*}
&& \bigl| R_\varepsilon(u_\varepsilon,u_\varepsilon) - R_\varepsilon
\bigl(u_\varepsilon, u_\varepsilon^{(4)}\bigr) \bigr|
\\
&&\qquad \leq C_K \int_\R\varepsilon^2 y^3
\bigl( |\hD_{\varepsilon y} u_\varepsilon|^2
+ |\hD_{\varepsilon y} u_\varepsilon| \bigl|\hD_{\varepsilon y}
u_\varepsilon^{(4)}\bigr|
+ \bigl|\hD_{\varepsilon y} u_\varepsilon^{(4)}\bigr|^2 \bigr)
\\
&&\qquad\quad\hspace*{26.6pt}{} \times\bigl|\hD_{\varepsilon y} \bigl(u_\varepsilon
-u_\varepsilon^{(4)} \bigr)\bigr| |\mu|(dy)
\\
&&\qquad \leq C_K \int_\R\varepsilon y^2
\bigl( |\hD_{\varepsilon y} u_\varepsilon| + \bigl|\hD_{\varepsilon y}
u_\varepsilon^{(4)}\bigr| \bigr)
\bigl|\hD_{\varepsilon y} \bigl(u_\varepsilon-u_\varepsilon^{(4)} \bigr)\bigr| |\mu
|(dy) .
\end{eqnarray*}
In view of (\ref{eqTheta-diff}) it thus follows that
\begin{eqnarray*}
&& \bigl\| R_\varepsilon(u_\varepsilon,u_\varepsilon) - R_\varepsilon
\bigl(u_\varepsilon, u_\varepsilon^{(4)}\bigr) \bigr\|_{L^1}
\\
&&\qquad \leq C_K \varepsilon\int_\R y^2
\bigl\|\bigl( |\hD_{\varepsilon y} u_\varepsilon| + \bigl|\hD_{\varepsilon y}
u_\varepsilon^{(4)}\bigr| \bigr)
\bigl| \hD_{\varepsilon y} \bigl(u_\varepsilon-u_\varepsilon^{(4)} \bigr)\bigr| \bigr\|
_{L^1} |\mu|(dy)
\\
&&\qquad \leq C_K \varepsilon\int_\R y^2
\bigl( \|\hD_{\varepsilon y} u_\varepsilon\|_{L^2}
+ \bigl\|\hD_{\varepsilon y} u_\varepsilon^{(4)}\bigr\|_{L^2} \bigr)
\bigl\| \hD_{\varepsilon y} \bigl(u_\varepsilon-u_\varepsilon^{(4)}\bigr) \bigr\|_{L^2}
|\mu|(dy)
\\
&&\qquad \leq C_K \varepsilon\sqrt{\bigl(\Theta_\varepsilon
(u_\varepsilon)
+\Theta_\varepsilon\bigl(u_\varepsilon^{(4)} \bigr)\bigr)\Theta_\varepsilon
\bigl(u_\varepsilon-u_\varepsilon^{(4)} \bigr)} .
\end{eqnarray*}
Using (\ref{eqTheta}), (\ref{eqTheta-rho}), and the definition of
the stopping time to bound $\Theta_\varepsilon(\psitlx)$, we find that
\begin{eqnarray*}
\Theta_\varepsilon(u_\varepsilon) + \Theta_\varepsilon
\bigl(u_\varepsilon^{(4)} \bigr)
& \leq & C\bigl( \Theta_\varepsilon( u_\varepsilon) + \Theta
_\varepsilon(\varrho_\varepsilon)
+ \Theta_\varepsilon(\psitlx)\bigr)
\\& \leq &C_K ( \varepsilon^{-1 - \kappa} + \varepsilon^{2\alpha-2}
+ \varepsilon^{\chi-1 - \kappa})
\leq\varepsilon^{-1 - \kappa} .
\end{eqnarray*}
Using (\ref{eqTheta-diff}), we thus obtain
%
%
\begin{equation} \label{eqone}
\bigl\| R_\varepsilon(u_\varepsilon,u_\varepsilon) - R_\varepsilon
\bigl(u_\varepsilon, u_\varepsilon^{(4)}\bigr) \bigr\|_{L^1}
\leq C_K ( \varepsilon^{\chi/2- 1/2 - 2\kappa}
+ \varepsilon^{\alpha-1/2 - \kappa} \|\varrho_\varepsilon\|
_\alpha) .\hspace*{-40pt}
\end{equation}
Furthermore, taking into account that
\[
\bigl\| u_\varepsilon- u_\varepsilon^{(4)}\bigr\|_{L^\infty}
\leq C_K ( \| \psitlx\|_{L^\infty}
+ \| \varrho_\varepsilon\|_{L^\infty})
\leq C_K ( \varepsilon^{\chi/2- \kappa}
+ \| \varrho_\varepsilon\|_{\alpha}) ,
\]
we have by (\ref{eqTheta}),
%
%
\begin{eqnarray}
\label{eqtwo}
&&
\bigl\| R_\varepsilon\bigl(u_\varepsilon,u_\varepsilon^{(4)}\bigr)
- R_\varepsilon\bigl(u_\varepsilon^{(4)},u_\varepsilon^{(4)}\bigr) \bigr\|_{L^1}
\nonumber\\
&&\qquad \leq C_K \varepsilon^2 \bigl\| u_\varepsilon- u_\varepsilon
^{(4)} \bigr\|
_{L^\infty}
\int_\R y^3
\bigl\| \bigl|\hD_{\varepsilon y} u_\varepsilon^{(4)}\bigr|^3 \bigr\|_{L^1}
|\mu|(dy)
\nonumber\\
&&\qquad \leq C_K \varepsilon\bigl\| u_\varepsilon- u_\varepsilon^{(4)} \bigr\|
_{L^\infty}
\int_\R y^2 \bigl\| \hD_{\varepsilon y} u_\varepsilon^{(4)} \bigr\|_{L^2}^2
|\mu|(dy)
\\
&&\qquad = C_K \varepsilon\bigl\| u_\varepsilon- u_\varepsilon^{(4)}
\bigr\|
_{L^\infty}\Theta_\varepsilon^2\bigl( u_\varepsilon^{(4)} \bigr)
\nonumber\\
&&\qquad \leq C_K ( \varepsilon^{\chi/2- 3\kappa}
+ \varepsilon^{- 2\kappa} \|
\varrho_\varepsilon\|_{\alpha}) .\nonumber
\end{eqnarray}
The claim follows by adding (\ref{eqone}) and (\ref{eqtwo}) and
using the embedding $L^1 \subseteq H^{-1}$.

\subsubsection*{Fifth term}
As in the first step, we have
%
%
\begin{eqnarray}
\label{eqv4-5}
\| \sigma_{\varepsilon,4} \|_{-1}
&=&
\| F(\vtg+ \psitg)
- F(\vtg+ \psitgx- \varrho_\varepsilon) \|_{-1}
\nonumber\\
& \leq &C_K ( \| \psitlx\|_{L^\infty}
+ \| \varrho_\varepsilon\|_{L^\infty} )
\\
& \leq &C_K ( \varepsilon^{\chi/2- \kappa}
+ \| \varrho_\varepsilon\|_{\alpha} ) .\nonumber
\end{eqnarray}

Combining the five estimates, we obtain
\begin{eqnarray*}
\|\varrho_\varepsilon(t)\|_{\alpha}
& \leq & C (t-s)^{(1-\alpha)/2}
\sup_{r\in(s,t)} \| \sigma_\varepsilon(r)\|_{-1}
\\& \leq & C (t-s)^{(1-\alpha)/2}
\sup_{r\in(s,t)}
\bigl( \varepsilon^{- 2 \kappa} \| \varrho_\varepsilon(r) \|
_{\alpha}
+ \varepsilon^{ \chi/2- 1/2 - 2 \kappa} \bigr) .
\end{eqnarray*}
The result now follows as in the proof of Proposition~\ref{propest-01}.
\end{pf}

\section*{Acknowledgments}

We are grateful to Hendrik Weber, Jochen Vo\ss\ and Andrew Stuart for
numerous discussions on this and related problems.
We thank the anonymous referee for useful comments.
Part of this work has been carried\vadjust{\goodbreak} out while J. Maas was visiting the
Courant Institute and the University of Warwick. He thanks both
institutions for their kind hospitality and support.



%
\printaddresses


\begin{thebibliography}{18}

\bibitem{AF03}
\begin{bbook}[mr]
\bauthor{\bsnm{Adams},~\bfnm{Robert~A.}\binits{R.~A.}} \AND
  \bauthor{\bsnm{Fournier},~\bfnm{John J.~F.}\binits{J.~J.~F.}}
(\byear{2003}).
\btitle{Sobolev Spaces},
\bedition{2nd} ed.
\bseries{Pure and Applied Mathematics (Amsterdam)}
\bvolume{140}.
\bpublisher{Elsevier}, \baddress{Amsterdam}.
\bid{mr={2424078}}
\end{bbook}
\endbibitem

\bibitem{BCJ94}
\begin{barticle}[mr]
\bauthor{\bsnm{Bertini},~\bfnm{L.}\binits{L.}},
  \bauthor{\bsnm{Cancrini},~\bfnm{N.}\binits{N.}} \AND
  \bauthor{\bsnm{Jona-Lasinio},~\bfnm{G.}\binits{G.}}
(\byear{1994}).
\btitle{The stochastic {B}urgers equation}.
\bjournal{Comm. Math. Phys.}
\bvolume{165}
\bpages{211--232}.
\bid{issn={0010-3616}, mr={1301846}}
\end{barticle}
\endbibitem

\bibitem{BertGia97}
\begin{barticle}[mr]
\bauthor{\bsnm{Bertini},~\bfnm{Lorenzo}\binits{L.}} \AND
  \bauthor{\bsnm{Giacomin},~\bfnm{Giambattista}\binits{G.}}
(\byear{1997}).
\btitle{Stochastic {B}urgers and {KPZ} equations from particle systems}.
\bjournal{Comm. Math. Phys.}
\bvolume{183}
\bpages{571--607}.
\bid{doi={10.1007/s002200050044}, issn={0010-3616}, mr={1462228}}
\end{barticle}
\endbibitem

\bibitem{BCF91}
\begin{barticle}[mr]
\bauthor{\bsnm{Brze{\'z}niak},~\bfnm{Z.}\binits{Z.}},
  \bauthor{\bsnm{Capi{\'n}ski},~\bfnm{M.}\binits{M.}} \AND
  \bauthor{\bsnm{Flandoli},~\bfnm{F.}\binits{F.}}
(\byear{1991}).
\btitle{Stochastic partial differential equations and turbulence}.
\bjournal{Math. Models Methods Appl. Sci.}
\bvolume{1}
\bpages{41--59}.
\bid{doi={10.1142/S0218202591000046}, issn={0218-2025}, mr={1105007}}
\end{barticle}
\endbibitem

\bibitem{Upwind}
\begin{barticle}[mr]
\bauthor{\bsnm{Courant},~\bfnm{Richard}\binits{R.}},
  \bauthor{\bsnm{Isaacson},~\bfnm{Eugene}\binits{E.}} \AND
  \bauthor{\bsnm{Rees},~\bfnm{Mina}\binits{M.}}
(\byear{1952}).
\btitle{On the solution of nonlinear hyperbolic differential equations by
  finite differences}.
\bjournal{Comm. Pure Appl. Math.}
\bvolume{5}
\bpages{243--255}.
\bid{issn={0010-3640}, mr={0053336}}
\end{barticle}
\endbibitem

\bibitem{DDT94}
\begin{barticle}[mr]
\bauthor{\bsnm{Da~Prato},~\bfnm{Giuseppe}\binits{G.}},
  \bauthor{\bsnm{Debussche},~\bfnm{Arnaud}\binits{A.}} \AND
  \bauthor{\bsnm{Temam},~\bfnm{Roger}\binits{R.}}
(\byear{1994}).
\btitle{Stochastic {B}urgers' equation}.
\bjournal{NoDEA Nonlinear Differential Equations Appl.}
\bvolume{1}
\bpages{389--402}.
\bid{doi={10.1007/BF01194987}, issn={1021-9722}, mr={1300149}}
\end{barticle}
\endbibitem

\bibitem{ZDP}
\begin{bbook}[mr]
\bauthor{\bsnm{Da~Prato},~\bfnm{Giuseppe}\binits{G.}} \AND
  \bauthor{\bsnm{Zabczyk},~\bfnm{Jerzy}\binits{J.}}
(\byear{1992}).
\btitle{Stochastic Equations in Infinite Dimensions}.
\bseries{Encyclopedia of Mathematics and Its Applications}
\bvolume{44}.
\bpublisher{Cambridge Univ. Press}, \baddress{Cambridge}.
\bid{doi={10.1017/CBO9780511666223}, mr={1207136}}
\end{bbook}
\endbibitem

\bibitem{DVS93}
\begin{barticle}[mr]
\bauthor{\bsnm{DeVore},~\bfnm{Ronald~A.}\binits{R.~A.}} \AND
  \bauthor{\bsnm{Sharpley},~\bfnm{Robert~C.}\binits{R.~C.}}
(\byear{1993}).
\btitle{Besov spaces on domains in {${\bf R}\sp d$}}.
\bjournal{Trans. Amer. Math. Soc.}
\bvolume{335}
\bpages{843--864}.
\bid{doi={10.2307/2154408}, issn={0002-9947}, mr={1152321}}
\end{barticle}
\endbibitem

\bibitem{Gyo98}
\begin{barticle}[mr]
\bauthor{\bsnm{Gy{\"o}ngy},~\bfnm{Istv{\'a}n}\binits{I.}}
(\byear{1998}).
\btitle{Existence and uniqueness results for semilinear stochastic partial
  differential equations}.
\bjournal{Stochastic Process. Appl.}
\bvolume{73}
\bpages{271--299}.
\bid{doi={10.1016/S0304-4149(97)00103-8}, issn={0304-4149}, mr={1608641}}
\end{barticle}
\endbibitem

\bibitem{Hai09}
\begin{bmisc}[auto:STB|2011-03-03|12:04:44]
\bauthor{\bsnm{Hairer},~\bfnm{M.}\binits{M.}}
(\byear{2009}).
\bhowpublished{An introduction to stochastic PDEs. Available at
  \href{http://arxiv.org/abs/arXiv:0907.4178}{arXiv:}
  \href{http://arxiv.org/abs/arXiv:0907.4178}{0907.4178}.}
\end{bmisc}
\endbibitem

\bibitem{Hai10}
\begin{bmisc}[auto:STB|2011-03-03|12:04:44]
\bauthor{\bsnm{Hairer},~\bfnm{M.}\binits{M.}}
(\byear{2011}).
\bhowpublished{Singular perturbations to semilinear stochastic heat equations.
  \textit{Probab. Theory Related Fields}. To appear. Available at
  \href{http://arxiv.org/abs/arXiv:1002.3722}{arXiv:1002.3722}.}
\end{bmisc}
\endbibitem

\bibitem{Jochen}
\begin{bmisc}[auto:STB|2011-03-03|12:04:44]
\bauthor{\bsnm{Hairer},~\bfnm{M.}\binits{M.}} \AND
  \bauthor{\bsnm{Voss},~\bfnm{J.}\binits{J.}}
(\byear{2010}).
\bhowpublished{Approximations to the stochastic {B}urgers equation.
\textit{J. Nonlinear Sci.} To appear.
Available at \href{http://arxiv.org/abs/arXiv:1005.4438}{arXiv:1005.4438}.}
\end{bmisc}
\endbibitem

\bibitem{KW92}
\begin{bbook}[mr]
\bauthor{\bsnm{Kwapie{\'n}},~\bfnm{Stanis{\l}aw}\binits{S.}} \AND
  \bauthor{\bsnm{Woyczy{\'n}ski},~\bfnm{Wojbor~A.}\binits{W.~A.}}
(\byear{1992}).
\btitle{Random Series and Stochastic Integrals: Single and Multiple}.
\bpublisher{Birkh\"auser}, \baddress{Boston, MA}.
\bid{mr={1167198}}
\end{bbook}
\endbibitem

\bibitem{Lunardi}
\begin{bbook}[mr]
\bauthor{\bsnm{Lunardi},~\bfnm{Alessandra}\binits{A.}}
(\byear{1995}).
\btitle{Analytic Semigroups and Optimal Regularity in Parabolic Problems}.
\bseries{Progress in Nonlinear Differential Equations and Their
Applications}
\bvolume{16}.
\bpublisher{Birkh\"auser}, \baddress{Basel}.
\bid{mr={1329547}}
\end{bbook}
\endbibitem

\bibitem{BookNA}
\begin{bbook}[mr]
\bauthor{\bsnm{Mattheij},~\bfnm{R.~M.~M.}\binits{R.~M.~M.}},
  \bauthor{\bsnm{Rienstra},~\bfnm{S.~W.}\binits{S.~W.}} \AND
  \bauthor{\bparticle{ten}
  \bsnm{Thije~Boonkkamp},~\bfnm{J.~H.~M.}\binits{J.~H.~M.}}
(\byear{2005}).
\btitle{Partial Differential Equations: Modeling, Analysis, Computation}.
\bpublisher{SIAM},
  \baddress{Philadelphia, PA}.
\bid{mr={2178809}}
\end{bbook}
\endbibitem

\bibitem{RevYor}
\begin{bbook}[mr]
\bauthor{\bsnm{Revuz},~\bfnm{Daniel}\binits{D.}} \AND
  \bauthor{\bsnm{Yor},~\bfnm{Marc}\binits{M.}}
(\byear{1994}).
\btitle{Continuous Martingales and {B}rownian Motion},
\bedition{2nd} ed.
\bseries{Grundlehren der Mathematischen Wissenschaften [Fundamental Principles
  of Mathematical Sciences]}
\bvolume{293}.
\bpublisher{Springer}, \baddress{Berlin}.
\bid{mr={1303781}}
\end{bbook}
\endbibitem

\bibitem{WZ}
\begin{barticle}[mr]
\bauthor{\bsnm{Wong},~\bfnm{Eugene}\binits{E.}} \AND
  \bauthor{\bsnm{Zakai},~\bfnm{Moshe}\binits{M.}}
(\byear{1965}).
\btitle{On the convergence of ordinary integrals to stochastic integrals}.
\bjournal{Ann. Math. Statist.}
\bvolume{36}
\bpages{1560--1564}.
\bid{issn={0003-4851}, mr={0195142}}
\end{barticle}
\endbibitem

\end{thebibliography}
\end{document}